\documentclass{article}[12pt]
\usepackage{amsmath,amsfonts,amsthm}
\usepackage{amssymb}

\usepackage{epsfig}

\usepackage{graphics}
\usepackage{graphicx}

\usepackage{latexsym}
\usepackage{graphics}
\usepackage{amsmath}
\usepackage{amsthm}
\usepackage{xspace}
\usepackage{amssymb}
\usepackage{color}
\usepackage{cite}
\usepackage{latexsym}
\usepackage{graphics}
\usepackage{amsmath}
\usepackage{amsthm}
\usepackage{xspace}
\usepackage{amssymb}
\usepackage{epsfig}

\newcommand{\beql}[1]{\begin{equation}\label{#1}}
\newcommand{\eeql}{\end{equation}}
\newcommand{\eqn}[1]{(\ref{#1})}

\newcommand{\R}{\mathbb{R}}
\newcommand{\pr}{\mathbb{P}}
\newcommand{\E}{\mathbb{E}}

\newcommand{\cs}{{\cal S}}

\newcommand{\cn}{{\cal N}}
\newcommand{\cj}{{\cal J}}

\newcommand{\calr}{{\cal R}}

\newcommand{\Z}{\mathbb{Z}}
\newcommand{\barZ}{\bar{\Z}}

\newtheorem{thm}{Theorem}
\newtheorem{lem}[thm]{Lemma}

\newtheorem{cor}[thm]{Corollary}

\newtheorem{definition}[thm]{Definition}

\newtheorem{rem}{Remark}

\begin{document}

\title{Pull-based load distribution 
among heterogeneous parallel servers: the case of multiple routers
}

\author
{
Alexander L. Stolyar \\
Lehigh University\\
200 West Packer Avenue, Room 484\\
Bethlehem, PA 18015 \\
\texttt{stolyar@lehigh.edu}
}

\date{\today}

\maketitle

\begin{abstract}

The model is a service system, consisting of several large server pools.
A server processing speed and buffer size (which may be finite or infinite) depend on the pool.
The input flow of customers is split equally among a fixed number of routers, which 
must assign customers to the servers immediately upon arrival.
We consider an asymptotic regime in which the customer total arrival rate
and pool sizes scale to infinity simultaneously, in proportion 
to a scaling parameter $n$, while the number of routers remains fixed.

We define and study a multi-router generalization of the pull-based customer assignment (routing) algorithm PULL, 
introduced in \cite{St2014_pull} for the single-router model.
Under PULL algorithm, when a server becomes idle it send a ``pull-message'' to a randomly uniformly selected router; each router operates independently -- it assigns an arriving customer
to a server according to a randomly uniformly chosen available (at this router) pull-message, if there is any,
or to a randomly uniformly selected server in the entire system, otherwise.

Under Markov assumptions (Poisson arrival process and independent exponentially distributed service requirements),
and under sub-critical system load,
we prove asymptotic optimality of PULL: as $n\to\infty$, 
the steady-state probability of an arriving customer
experiencing blocking or waiting, vanishes.
Furthermore, PULL has an extremely low router-server message exchange rate of one message per customer.
These results generalize some of the single-router results in \cite{St2014_pull}.

\end{abstract}




\section{Introduction}

Design of efficient load distribution algorithm for modern data processing systems is a challenging problem.
An algorithm needs to have a good performance, which typically means low waiting times and blocking probabilities of the jobs. 
A big part of the challenge of achieving this goal efficiently,
stems from the fact that such systems are heterogeneous (i.e., different servers may have different capabilities) and have very large scale. This makes the following features of a load distribution algorithm very desirable, and perhaps even required:

\begin{itemize}
\item 
An algorithm should be oblivious of the server types as much as possible.
\item An algorithm should allow a distributed implementation, where assignment (routing) of jobs to servers is done by multiple routers, each handling a fraction of demand. This is because having a single router to handle a massive demand
may be infeasible or impractical (see \cite{G11}).
\item The router-server signaling overhead should be small.
\end{itemize}

Motivated by the challenges described above, in this paper we consider the following model.
It is a service system, consisting of several large server pools.
The system is heterogeneous in that
a server processing speed and buffer size
depend on the pool; the buffer size is the maximum number of jobs  -- or, customers -- that can be queued at the server, and it can be finite or infinite.
There is a finite number 
of {\em routers}. The flow of customer arrivals into the system is split equally among the routers, which 
must assign (or, route) ``their'' customers to the servers immediately upon arrival.
This model is a generalization of the single-router model in \cite{St2014_pull}.

We define and study (two versions of) a {\em pull-based} customer assignment (routing) algorithm PULL.
(One of the algorithm versions is a generalization of the single-router PULL algorithm in \cite{St2014_pull}.)
Under this algorithm, when a server becomes idle it sends a ``pull-message'' to a randomly uniformly selected router; each router operates independently -- it assigns an arriving customer
to a server according to a randomly uniformly chosen available (at this router) pull-message, if there is any,
or to a randomly uniformly selected server in the entire system, otherwise.

We make Markov assumptions: Poisson arrival process and independent exponentially distributed customer service requirements. We consider an asymptotic regime in which 
the number of routers remains fixed, but 
the customer total arrival rate
and pool sizes scale to infinity simultaneously, in proportion 
to some scaling parameter $n$, while the system load is
sub-critical (i.e., the load is smaller than the system capacity by $cn$ for some $c>0$).

The {\bf main results} of this paper can be informally stated as follows (formal statements are
given in Lemma~\ref{lem-stabil-finite}, Theorems~\ref{thm-stabil-infinite} and \ref{thm-infinite}):
\begin{itemize}
\item The system process under PULL is stable (positive recurrent) for all sufficiently large $n$.
\item PULL is asymptotically optimal: As $n\to\infty$, the steady-state probability of each pool having idle servers and each router having pull-messages goes to $1$; consequently,
the steady-state probability of an arriving customer
experiencing blocking or waiting, vanishes. (Furthermore, in the limit, the idle servers in each pool have their pull-messages distributed equally among the routers.)
\item Router-server signaling message exchange rate under PULL is at most two messages per customer for 
a pre-limit system, and equal to one message per customer in the $n\to\infty$ limit.
\end{itemize}

These results generalize the single-router results in \cite{St2014_pull}, under Markov assumptions.
(The main results in \cite{St2014_pull} hold under more general assumptions on the service requirement distributions -- it suffices that they have a {\em decreasing hazard rate}.)

The following are some comments on the model and main results, which demonstrate that PULL algorithm does satisfy the listed above desired features of a load distribution algorithm:
\begin{itemize}
\item The assumption that the input flow of customers is equally split among the routers is non-restrictive. A real system can have a single ``pre-router'', which is a single point of entry of customers into the system.
The pre-router's only function would be to distribute arriving customers in, say, round-robin fashion among the 
routers. Given the utmost simplicity of this function, there is no problem in handling all arriving customers
by a single pre-router;
this is in contrast to a router, which has to do actual customer-to-server assignment, requiring more
processing per customer.
\item The asymptotic optimality property of PULL shows that its performance in a large-scale system is perfect --
waiting and/or blocking of arriving customers vanishes.
In particular, its performance is much superior to that of the celebrated {\em power-of-d-choices}, or JSQ(d), algorithm \cite{VDK96,Mitz2001,BLP2012-jsq-asymp-indep,BLP2013-jsq-asymp-tail}. 
(JSQ(d) is defined and discussed below. See also \cite{St2014_pull} for more detailed PULL Vs. JSQ(d) comparisons.)
The advantage of PULL over JSQ(d) is especially striking in heterogeneous systems,
where JSQ(d) is not even applicable in general (see  \cite{St2014_pull}).
\item The router-server message exchange rate under PULL is extremely low: at most 2/customer and, in the limit, 1/customer. (It is $2d$ under JSQ(d), where $d$ is a parameter whose values of interest are $d\ge 2$.)
\item We want to emphasize that PULL algorithm does {\em not} try to explicitly balance the distribution of available pull-messages among routers -- a perfect (in the limit) balancing occurs ``automatically''. This useful property is somewhat counterintuitive.
\item In terms of technique, the generalization of the results of \cite{St2014_pull} to the multi-router case is {\em not straightforward}. The main reason is as follows. Process monotonicity plays a key role in this paper as it did in \cite{St2014_pull}. However, in the single-router case, there is the minimum (smallest) state of the process; this gives a simple characterization of the process steady-state as a limit of the stochastically
monotone increasing process, starting from the smallest state; the analysis in  \cite{St2014_pull} essentially relies
on this fact. There is {\em no} smallest state in the multi-router case; this requires new approaches to establish 
the main results.
\end{itemize}

There exists a large amount of previous work on the load distribution in service systems.
 Majority of the previous work is focused on
{\em load balancing} (in the sense of equalizing server loads)
in homogeneous systems (with all servers identical); see e.g. \cite{BLP2013-jsq-asymp-tail,G11} for reviews. Load balancing is only one possible objective of load distribution.
While it is a natural objective for homogeneous systems, 
it is not necessary to achieve asymptotic optimality.
Our PULL algorithm is an example -- it is asymptotically optimal, 
without attempting and without in general achieving load balancing
among all servers in the system. (It does, however, equalize load among servers 
within each pool.)

Much attention, especially in recent years, was devoted to {\em power-of-d-choices}, or {\em join-the-shortest-queue(d)}, or JSQ(d), algorithm \cite{VDK96,Mitz2001,BLP2012-jsq-asymp-indep,BLP2013-jsq-asymp-tail}.
Under this algorithm each arriving customer joins the shortest out of $d$ randomly selected queues,
where $d\ge 1$ is an algorithm parameter. JSQ(1) is just a random uniform routing, so that the interesting cases are when $d\ge 2$ and is small.
The cited above work, as well as practically all work on JSQ(d),
is for homogeneous systems. (For heterogeneous systems JSQ(d) with $d\ge 2$ is not applicable 
in general -- see \cite{St2014_pull}.) The advantage of JSQ(d) with small $d\ge 2$ over JSQ(1) is that it dramatically reduces the steady-state delays, at the cost of only small router-server signaling message exchange rate of $2d$/customer. Our results show that PULL algorithm asymptotically eliminates
steady-state delays and has even smaller signaling message exchange rate, just $1$/customer;
and it is applicable to heterogeneous systems. Assuming the queue-length queries/responses are instantaneous,
the JSQ(d) algorithm is unaffected by whether
the system has a single or multiple routers.

A pull-based approach to load distribution is relatively recent \cite{BB08,G11,St2014_pull}.
Paper \cite{BB08} proposed this approach in practical settings, and studied it by simulations,
which demonstrated good performance.
Paper \cite{G11} considers a homogeneous system model with multiple routers.
It proposes a pull-based algorithm, called {\em join-idle-queue} (JIQ),
which is essentially what we call PULL algorithm in this paper.
(PULL has additional features; most importantly, the random uniform choice of pull-messages by each router.)
The asymptotic regime in \cite{G11} is different from ours: both the numbers of servers and routers
grow to infinity, with their ratio kept constant; heuristic arguments are used to conjecture the limit 
of the system steady-states; simulations both supported the conjecture validity and demonstrated good performance of the algorithm. Paper \cite{St2014_pull} considered a special -- single-router -- version 
of the model in this paper, and rigorously proved the asymptotic optimality of PULL,
under assumptions more general than Markov (namely, for service time distributions with decreasing hazard rate).
The main results of this paper generalize those in \cite{St2014_pull} to the multi-router model, under Markov
assumptions.

Recent papers \cite{EsGam2015,MBLW2015} consider a single-router homogeneous system, under Markov
assumptions, in the so called Halfin-Whitt asymptotic regime, where the number of servers scales in proportion
to $n\to\infty$ and the system load is smaller than capacity by $c\sqrt{n}$ for some $c>0$.
These papers prove diffusion limits of the system transient behavior under {\em join-the-shortest-queue}
routing (in \cite{EsGam2015}) and under its generalization which also includes JIQ as a special case
(in \cite{MBLW2015}). 

Finally, we note that pull-based algorithms can be applied to systems with more general server structures, beyond simple single server. For example, it is shown in \cite{St2014_pull} that the asymptotic 
optimality of PULL extends to (single-router) models, where server processing speed depends on its queue length.
Another example is recent paper \cite{St2015_grand-het}, which proposes and studies 
an algorithm, which can be viewed as a version of PULL, for (single-router) heterogeneous service systems
with {\em packing constraints} at the servers.

\subsection{Basic notation}

Symbols $\R, \R_+, \Z, \Z_+$ denote the sets of real, real non-negative,
integer, and integer non-negative numbers, respectively.
For finite- or infinite-dimensional vectors, the vector inequalities 
are understood component-wise, unless explicitly stated otherwise.
We use notation $x(\cdot)=(x(t), ~t\ge 0)$ for both a random process
and its realizations, the meaning is determined by the context;
the state space (of a random process) and the metric and/or 
topology on it are defined where appropriate,
and we always consider Borel $\sigma$-algebra on the state space.
Notations $\Rightarrow$ and $\stackrel{d}{=}$
signify convergence and equality {\em in distribution}, respectively,
for random elements; $\stackrel{P}{\to}$ denotes convergence {\em in probability} of random variables.
For a process $x(\cdot)$, we denote by $x(\infty)$ a random element
whose distribution is a {\em stationary distribution} of the process,
assuming the latter exists. 
For $a,b\in \R$, $\lfloor a \rfloor$ denotes the largest integer less than 
or equal to $a$, $a\wedge b = \min\{a,b\}$.
We use notation $Dist[X]$ for the probability distribution
of a random element $X$; for a measure $\chi$ its total variation is denoted $\|\chi\|$.
Left limit of a function (or a process) $z(t)$ at a given point $t$ is denoted $z(t-)$;
the right derivative of a function at $t$ is denoted  $(d^+/dt) z(t)$. 
Abbreviation {\em u.o.c.} means 
{\em uniform on compact sets} convergence; {\em w.p.1} means
{\em with probability 1}; {\em a.e.} means 
{\em almost everywhere} with respect to Lebesgue measure;
 {\em WLOG} means 
{\em without loss of generality}; {\em RHS} and {\em LHS} mean {\em right-hand side}
and {\em left-hand side}, respectively.

\subsection{Paper organization} 

The model and main results on stability (Lemma~\ref{lem-stabil-finite} and Theorem~\ref{thm-stabil-infinite})
and asymptotic optimality (Theorem~\ref{thm-infinite}) of PULL are given  in Section~\ref{sec-model},
which also contains a discussion of PULL implementation mechanisms.
Section~\ref{sec-monotonicity} introduces an order relation on the
process state space, and establishes monotonicity properties, which are a key tool for our analysis;
it also discusses differences in implications of monotonicity in the single-router and multi-router cases.
In Section~\ref{sec-fluid-many-server} we study the process fluid limits, which are another important tool used in the paper.
We study a special system where all buffer sizes are $1$ in Section~\ref{sec-unit-buffers};
these are auxiliary results, needed for analysis of the general case.
Proofs of Theorems~\ref{thm-stabil-infinite} and \ref{thm-infinite} are in Sections~\ref{sec-stability}
and \ref{sec-proof-finite-buffers}, respectively.
We conclude in Section~\ref{sec-conclusion}.

\section{Model and main results}
\label{sec-model}

\subsection{Model structure}
\label{subsec-model}

{\bf Arrival process.} Customers  arrive according to a Poisson 
process of rate $\Lambda>0$.

{\bf Servers.} 
There are $J\ge 1$ server pools. Pool $j \in \cj \equiv \{1,\ldots,J\}$
consists of $N_j$ identical servers.
Servers in pool 1 are indexed by $i \in \cn_1 =\{1,\ldots, N_1\}$,
in pool 2 by $i \in \cn_2 =\{N_1+1,\ldots, N_1+N_2\}$, and so on;
$\cn = \cup \cn_j$ is the set of all servers.
The service time of a customer at a server in pool $j$
is an independent, exponentially distributed random variable 
 with mean $1/\mu_j\in (0,\infty)$, $j \in \cj$.
We assume that the customers at any server are served in the 
first-come-first-serve (FCFS) order. 
(That is, at any time
only the head-of-the-line customer at each server is served.)
The buffer size (maximum queue length) at any server in pool $j$ is $B_j\ge 1$;
we allow the buffer size to be either finite, $B_j < \infty$, or infinite, $B_j=\infty$.
A new customer, routed to a server $i\in \cn_j$, joins the queue at that server,
unless $B_j$ is finite and the queue length $Q_i=B_j$ -- in this case 
the customer is blocked
(or lost; i.e., it leaves the system immediately, without receiving any service).

{\bf Routing.} We assume that the customer arrival process is distributed equally among $R$ routers,
labeled by $r\in \calr = \{1,\ldots,R\}$. 
Specifically, each arriving customer is assigned randomly uniformly to one of the routers; therefore, the arrival processes to the routers are
independent Poisson processes of rate $\Lambda/R$. When a router receives an arriving customer, 
it has to either immediately route it to one of the servers or immediately block it from service (in which case the customer is lost).

\subsection{Asymptotic regime}
\label{sec-asymptotic-reg}

We consider the following (many-servers) asymptotic regime. 
The total number of servers $n=\sum_j N_j$ is the scaling parameter, which
increases to infinity; the arrival rate and the server pool sizes increase in proportion 
to $n$, $\Lambda= \lambda n$, $N_j=\beta_j n, ~j\in \cj$, where 
$\lambda, \beta_j, j\in \cj$, are positive constants, $\sum_j \beta_j = 1$.
(To be precise, the values of $N_j$ need to be integer, 
e.g. $N_j=\lfloor \beta_j n \rfloor$. Such definition
would not cause any problems, besides clogging
notation, so we will simply assume that all 
$\beta_j n$ ``happen to be'' integer.)
We assume that the subcritical load condition holds:
\beql{eq-load}
\lambda < \sum_j \beta_j \mu_j.
\end{equation}

\subsection{PULL routing algorithm}
\label{sec-alg}

We study the following pull-based algorithm, labeled PULL.

\begin{definition}[PULL]
\label{def-pull-basic1}
We consider two versions of the PULL algorithm, labeled PULL-1 and PULL-2.
They have common components (a)-(d) below; the components (e.1) and (e.2) pertain only to 
PULL-1 and PULL-2, respectively.
\begin{itemize}
\item[(a)] At any given time,
each idle server has exactly one {\em pull-message} (containing the server identity), 
located at one of the routers. (Equivalently, each idle server is {\em associated} with one of the routers,
and each router ``knows'' the subset of idle servers currently associated with it.) A pull-message located
at router $r$ is called $r$-pull-message.
\item[(b)] When a server becomes idle after a service completion and departure of a customer, a pull-message from this server is generated and placed at one of the routers, chosen randomly uniformly.
\item[(c)] When a customer arrives at an idle server, the pull-message for this server is destroyed.
(The pull-message is destroyed whether it was located at the router which sent the customer or at a different router.)
\item[(d)] When a customer arrives at router $r$, if there are available $r$-pull-messages, 
the customer is sent to one of the corresponding idle servers, chosen randomly uniformly.
\item[(e.1)] Under PULL-1, when a customer arrives at router $r$, if there are no available $r$-pull-messages, 
the customer is blocked (and lost).
\item[(e.2)] Under PULL-2, when a customer arrives at a router, if there are no available $r$-pull-messages, 
the customer is 
sent to one of the servers in the entire system, chosen randomly uniformly. (Recall that, according to the 
model structure, the customer then may be blocked at the server; this occurs if and only if the server's queue is ``full'', i.e. the queue length is equal to the finite buffer size.)
\end{itemize}
\end{definition}

Note that, if we are interested in the system stationary distributions, then
 for PULL-1, WLOG, we can assume $B_j=1$ for all $j$, because for any actual buffer sizes
and any initial state, w.p.1 after a finite time there will be at most one customer at each server. We adopt this assumption for PULL-1 for the rest of the paper. Consequently, for any system for which condition 
$B_j=1, \forall j$, does {\em not} hold, we consider PULL-2 algorithm only.

We also remark that the PULL algorithm in \cite{St2014_pull} is a special case (for $R=1$) of PULL-2, but not of PULL-1.

\subsection{Main results}
\label{sec-result}

Our main results hold for both versions, PULL-1 and PULL-2, of the algorithm. (Of course, this does {\em not} mean
that the system {\em behavior} is same under both versions.) In general, throughout the paper,
when we say that a certain fact holds for PULL algorithm, we mean that it holds for any fixed version of it,
i.e. PULL-1 or PULL-2.

In the system with parameter $n$, the system state is given by 
$S^n=(Q^n_i, ~ D^n_i, ~i\in \cn)$, where the components are defined as follows:
$Q^n_i \in \Z_+$ is the
queue length at server $i$; $D^n_i \in \calr$ 
is the label of the router 
that ``holds'' the pull-message from server $i$, if it is idle ($Q^n_i =0$), and $D^n_i = 0$, otherwise.
Obviously, any state $S^n$ satisfies condition $I\{Q^n_i=0\} = I\{D^n_i \ne 0\}, ~~\forall i\in \cn$.

Due to symmetry of servers within each pool,
the alternative -- {\em mean field}, or {\em fluid-scale}  
-- representation of the process is
as follows. 
Let $x^n_{k,j}$ denote the fraction of the (total number of) servers,
which are in pool $j$ and have queue length greater than or equal to $k$;
let $\xi^n_{r,j}$ denote the fraction of the (total number of) servers,
which are idle in pool $j$ and have a pull-message at router $r$.
Denote 
$$
x^n = (x^n_{k,j}, ~k\in \Z_+, ~j\in \cj),~~
\xi^n = (\xi^n_{r,j}, ~r\in \calr, ~j\in \cj).
$$
Then 
$$
s^n = (x^n, \xi^n),
$$
is (the mean field representation of ) the system state.
We will view states
$s^n$, for any $n$, as elements of the common space
with elements denoted by 
$$
s = (x, \xi), ~~~~ x = (x_{k,j}, ~k\in \Z_+, ~j\in \cj),~~
\xi = (\xi_{r,j}, ~r\in \calr, ~j\in \cj).
$$
Specifically, the common state space is
$$
\cs = \{s~|~ 
\beta_j = x_{0,j} \ge x_{1,j} \ge x_{2,j} \ge \cdots \ge 0 ~\mbox{and}~ x_{0,j} - x_{1,j} = \sum_r \xi_{r,j}, 
\forall j\in\cj\},
$$
equipped with metric
\beql{eq-metric}
\rho(s,s') = \sum_{j} \sum_k 2^{-k} \frac{|x_{k,j}-x'_{k,j}|}{1+|x_{k,j}-x'_{k,j}|}
+ \sum_{j} \sum_r  |\xi_{r,j}-\xi'_{r,j}|,
\end{equation}
and the corresponding
Borel $\sigma$-algebra. (The corresponding topology is that of componentwise convergence.) 
Space $\cs$ is compact. 

For any $n$, the process $S^n(t), ~t\ge 0,$ -- and its projection $s^n(t), ~t\ge 0,$ -- 
is a continuous-time, 
countable state space, irreducible Markov process. 
(For any $n$, the state space of $s^n(\cdot)$ is a countable subset of $\cs$.)
If the buffer sizes $B_j$ are finite in {\em all} pools $j$, 
the state space is obviously finite, and therefore the process $S^n(\cdot)$
(and then $s^n(\cdot)$)
is positive recurrent (and then ergodic), with unique stationary distribution. 
Therefore, the following fact is trivially true -- we present it for 
ease of reference.

\begin{lem}
\label{lem-stabil-finite}
Consider the system under PULL algorithm. If $B_j<\infty$ for all $j\in\cj$, the Markov process $s^n(\cdot)$ is ergodic.
\end{lem}

Our first main result is that $s^n(\cdot)$ is ergodic under arbitrary (possibly infinite) buffer sizes,
when $n$ is large enough.

\begin{thm}
\label{thm-stabil-infinite}
Consider the system under PULL algorithm. Suppose that $B_j =\infty$ for some (or all) $j\in\cj$.
(Recall that in this case we only need to consider PULL-2, because under PULL-1 all $B_j=1$, WLOG.)
Then, for all sufficiently large $n$, the Markov process $s^n(\cdot)$ is ergodic.
\end{thm}

Define numbers $\nu_j \in (0,\beta_j)$, $j\in \cj$,
uniquely determined by the conditions
\beql{eq-invar-point1}
\lambda = \sum_j \nu_j \mu_j, 
~~~ \nu_j \mu_j / (\beta_j - \nu_j) = \nu_\ell \mu_\ell / (\beta_\ell - \nu_\ell),
~\forall j,\ell\in\cj.
\end{equation}
Let us define the {\em equilibrium point} $s^* \in \cs$ by
\beql{eq-invar-point2}
x^*_{1,j} = \nu_j, ~~ x^*_{k,j} = 0,~k \ge 2, ~~~ j\in\cj,
\end{equation}
\beql{eq-invar-point3}
\xi^*_{r,j} = (\beta_j - \nu_j)/R, ~~~~ r \in \calr, ~ j\in\cj.
\end{equation}

The equilibrium point $s^*$ definition in \eqn{eq-invar-point1}-\eqn{eq-invar-point3}
has the following interpretation. Condition \eqn{eq-invar-point2} means that point $s^*$ is such that the fraction $\nu_j < \beta_j$ of servers
(out of the total number of servers) in pool $j$ is occupied by exactly one customer, while the remaining servers in pool $j$ are idle.
Moreover, according to \eqn{eq-invar-point3},
{\em the idle servers is each pool are in equal quantities associated with (have a pull-message at) different routers;} this in particular means that each router has pull-messages available, and therefore all new arrivals are
routed according to pull-messages. Finally,
the numbers $\nu_j$ are uniquely determined by the condition \eqn{eq-invar-point1}, which simply says that the rate at which 
new arrivals are routed to pool $j$ (it is proportional to $\beta_j - \nu_j$) is equal to the service completion rate at pool $j$ 
(it is proportional to $\nu_j \mu_j$).

By Lemma~\ref{lem-stabil-finite} and Theorem~\ref{thm-stabil-infinite} process $s^n(\cdot)$
is ergodic for all large $n$ (and if all buffers are finite -- for all $n$).
Denote by $s^n(\infty)$ a random element with the distribution
equal to the stationary distribution of $s^n(\cdot)$.

Our second main result is the following 
\begin{thm}
\label{thm-infinite}
Under PULL algorithm,
$s^n(\infty)\Rightarrow s^*$ as $n\to\infty$.
\end{thm}

Given the definition of $s^*$, Theorem~\ref{thm-infinite} implies that,
in the $n\to\infty$ limit, there are always pull-messages available at each router.
(Furthermore, the pull-messages from idle servers in each pools are distributed equally among the routers.)
Therefore, PULL algorithm is {\em asymptotically optimal} in the following sense:
as $n\to\infty$, the steady-state probability of an arriving customer
experiencing blocking or waiting, vanishes. (Together, 
Theorems~\ref{thm-stabil-infinite} and \ref{thm-infinite}
generalize Theorem 2 in \cite{St2014_pull}  to the model with multiple routers.)
We discuss the implications of Theorem~\ref{thm-infinite} in more detail in Section~\ref{sec-practical}.

\subsection{PULL algorithm implementation mechanisms}
\label{sec-practical}

We will now discuss some natural {\em mechanisms} for the PULL algorithm practical implementation.
This discussion will demonstrate that very simple and extremely efficient mechanisms do exist, which is a 
major motivation for the PULL algorithm in the first place. It is also important to emphasize
that {\em as far as the algorithm analysis is concerned}, the actual implementation mechanism is, of course,
irrelevant -- Definition~\ref{def-pull-basic1} of the algorithm is sufficient for the analysis.

{\bf PULL-1.} When a server is initialized in the actual system, it is idle, and it sends a pull-message
to one of the routers chosen uniformly at random. After that,
the server sends a new pull-message to a randomly uniformly chosen router, 
after any service completion that leaves the server idle. (By definition of PULL-1, any service completion
leaves the server idle.) Each pull-message contains the server identity.
When a customer arrives at a router, if this router has available pull-messages, it 
picks one of those pull-messages
uniformly at random,  sends the customer to the corresponding server,
and immediately destroys the ``used'' pull-message.
When a customer arrives at a router, if this router has no available pull-messages,
the customer is blocked (and discarded).

{\bf PULL-2.} When a server is initialized in the actual system, it is idle, and it sends a pull-message
to one of the routers chosen uniformly at random. 
After that, the server sends a new pull-message to a randomly uniformly chosen router, 
after any service completion that leaves the server idle. (Under PULL-2, in general, not every service 
completion leaves the server idle.) 
Each pull-message contains the server identity.
Each time a server sends a pull-message, it ``remembers'' the identity
of the router, to which it was sent.
When a customer arrives at a router, if this router has available pull-messages, it 
picks one of those pull-messages
uniformly at random,  sends the customer to the corresponding server,
and immediately destroys the ``used'' pull-message.  
When a customer arrives at a router, if this router has no available pull-messages,
the customer is sent to one of the servers in the entire system, chosen uniformly at random.
Every customer sent from a router to a server is supplied 
with one additional bit of information, indicating whether or not a pull-message was used for its routing or not.
When a customer arrives to an idle server, the server checks whether or not a pull-message was used in routing;
if not, the server sends a special ``pull-remove-message'' (containing identity of the server) to the router which 
``held'' the pull-message from the server. When a router receives a pull-remove-message, it destroys the
pull-message of the corresponding server.

Next we discuss key properties of these implementation mechanisms.

\subsubsection{Asymptotic optimality is achieved by a very simple mechanism.} 

PULL algorithm implementation is clearly very simple. Routers do not keep track of the exact states of the servers. Each router only needs to have the list of servers (which is necessary
under any algorithm, as discussed in Section 2.5.2 of \cite{St2014_pull}), where for each server one additional bit is used to indicate the presence of a pull-message. Each server only needs to have the list of routers
(which it needs under any feedback-based algorithm), and one additional variable to ``remember'' which 
router currently holds its pull-message. 

Also note that there is {\em no} additional mechanism to balance the numbers of pull-messages
at different routers; such balancing (in the limit) occurs automatically. (A pull-based algorithm proposed 
in \cite{G11} for the case of multiple routers, does suggest a pull-message balancing mechanism. Our results
prove that when the ratio of servers to routers is large -- which is our asymptotic regime -- such additional mechanism is unnecessary.)

\subsubsection{Performance: message exchange rate.} 

It is clear that, under PULL-1, the steady-state rate of message exchange between the routers and the servers 
is {\em at most one} message per customer. Specifically, one pull-message is sent per each customer upon its service completion. (Some fraction of customers is blocked, so for a finite $n$ not every customer generates 
a pull-message.) As $n\to\infty$, the steady-state blocking probability vanishes; therefore, in the $n\to\infty$ limit, the message exchange rate is exactly one message per customer.

Under PULL-2, the steady-state router/server message exchange rate 
is {\em at most two} messages per customer. Specifically, 
for each customer, 
at most one pull-remove-message is sent upon its arrival at a server and 
at most one pull-message is sent upon its service completion (when it leaves the server idle).
In the $n\to\infty$ limit, the message exchange rate is only one message per customer (same as for PULL-1),
because the probability that a customer is routed {\em without} using a pull-message vanishes.

The asymptotic message exchange rate of one per customer is much lower than that
under JSQ(d) algorithm, for which it is $2d$. (See \cite{St2014_pull} for a more detailed discussion.)

\subsubsection{Performance: absence of routing delay.} 
\label{sec-routing-delay}

Pull-messages do {\em not} contribute
to the routing delay: an arriving customer is not waiting at the router for any pull-message or pull-remove-message, the routing decision is made 
immediately. This is unlike JSQ(d) algorithm, where each arriving customer waits for the queue-length request/response message exchange to complete, before being routed. (See also \cite{G11} for a discussion of this issue.)

\subsubsection{The notion of servers pools is purely logical.}

Neither routers nor servers use the notion of a server pool.
Indeed, each router need {\em not} know anything about each server parameters (service speed, buffer size) 
or state (queue length),
besides currently having or not a pull-message from it.
This means that
from the ``point of view'' 
of each router all servers form a single pool.
The notion of a server pool that we have in the model is purely logical, used for analysis only.
A real system may consist of a single or multiple pools of {\em non-identical} 
servers. In this case, we consider all servers of a particular type as forming a {\em logical} pool.
Our results still apply, as long as the number of servers of each type in the entire system is large.

\subsection{Algorithm implementation and physical system constraints/requirements}
\label{sec-practical2}

Clearly, the architecture and/or requirements of an actual physical system may impose constraints on 
any algorithm implementation. For example, if in the actual system ``servers'' simply cannot generate own messages to the routers and can only respond to the queue length queries from the routers, then PULL algorithm implementation mechanisms described above would not be feasible, while JSQ(d) mechanism could be.
 (This is typically not the case for cloud-based systems, as well as other web-based services, where servers have sufficient intelligence to generate own signaling messages; see, e.g. \cite{BB08,G11}.)
As another example, the routers may not have enough buffer space for the arriving customers to tolerate any routing delay (see Section~\ref{sec-routing-delay}) -- in this case JSQ(d) implementation would not be feasible, 
while PULL implementation could be; moreover, if routers have very limited intelligence and/or processing power, then PULL-1 might be the only option in terms of implementation. The applicability of JSQ(d) in heterogeneous systems
could be limited, because it does not ensure system stability (see \cite{St2014_pull}).
We believe that in most modern data processing systems the implementation of PULL is feasible,
but it is up to a practitioner to take into account all system constraints and choose the best algorithm 
for a specific application.

\section{The process in the compactified state space.
Monotonicity}
\label{sec-monotonicity}

All results in this section concern a system with a fixed $n$. Also, they hold for any fixed version of PULL -- either PULL-1 or PULL-2.

\subsection{Compactified state space. Order relation.}
\label{sec-compactified}

It will be convenient to consider a more general system and the Markov process.
Namely, we assume that the queue length at any server $i\in\cn_j$
within a pool $j$ with infinite buffer size ($B_j=\infty$), can be 
infinite. In other words, $Q_i^n(t)$ can take values in the
set $\barZ_+ \doteq \Z_+ \cup \{\infty\}$, which is the
one-point compactification of $\Z_+$, containing the ``point at infinity.''
We consider the natural topology and order relation on $\barZ_+$.
Obviously, $\barZ_+$ is compact. (Note that if $A$ is a finite subset
of $\Z_+$, then sets $A$ and $\barZ_+ \setminus A$ are both closed and open.)

Therefore, the state space of the generalized version of 
Markov process $S^n(\cdot)$ is the compact set $\barZ_+^n \times [\calr \cup \{0\}]^n$. 
The process transitions are defined in
exactly same way as before, with the following additional convention.
If $Q^n_i(t)=\infty$, then neither new arrivals into this queue nor
service completions in it, change the system state;
this means, in particular, that $Q^n_i(\tau)\equiv \infty$ for all $\tau\ge t$.

The corresponding generalized version of 
the process $s^n(\cdot)$ is defined as before; if at time $t$ some of the
queues in pool $j$ are infinite, then $x^n(t)$ is such that
$\lim_{k\to\infty} x^n_{k,j}(t) >0$. Note that the 
state space of the generalized $s^n(\cdot)$ is still the
compact set $\cs$, as defined above.

It is easy to see that, for each $n$, the (generalized versions of) 
processes $S^n(\cdot)$ and $s^n(\cdot)$ are Feller continuous.

We will consider the following natural order relation on the process $S^n(\cdot)$ state space.
Vector inequality $Q'\le Q''$ for $Q', Q'' \in \barZ_+^n$
is understood component-wise. If $D', D'' \in [\calr \cup \{0\}]^n$, then $D' \le D''$ is understood as the following
property:
\beql{eq-D-order}
D''_i \ne 0 ~~\mbox{implies}~~ D'_i = D''_i, ~~~\forall i \in \cn.
\end{equation}
Then, $S' = (Q',D') \le S''=(Q'',D'')$ is defined as $Q'\le Q''$ and $D'\le D''$.
The meaning of this order relation is clear: the state $S'$ is dominated by $S''$ if 
its queue length is smaller or equal at each server, and  the set of $r$-pull-messages in $S'$ is a superset
of that in $S''$ for each $r$ at each pool $j$.

The stochastic order relation $S'\le_{st} S''$ between two {\em random elements} 
means that they can be constructed on the
same probability space so that $S'\le S''$ holds w.p.1.

The corresponding order relation on the state space $\cs$ for the mean-field process $s^n(\cdot)$
is as follows. For $s'=(x',\xi'), s''=(x'',\xi'') \in \cs$, $s' \le s''$ is defined as $x' \le x''$ component-wise and
the property:
\beql{eq-xi-order}
\xi'_{r,j} \ge  \xi''_{r,j} ~~~\forall r \in \calr, ~\forall j\in \cj.
\end{equation}
The stochastic order relation $s' \le_{st} s''$ between two {\em random elements} in $\cs$
means that they can be constructed on the
same probability space so that $s'\le s''$ holds w.p.1.

In the rest of the paper, for a state $s^n(t)$ 
(with either finite $t\ge 0$ or $t=\infty$), we denote by
\beql{eq-x-inf}
x^n_{\infty,j}(t)\doteq \lim_{k\to\infty} x^n_{k,j}(t)
\end{equation}
the fraction of queues that are in pool $j$ and are infinite.
(Note that $x^n_{\infty,j}(t)$ is a function, but {\em not a component},
of $x^n(t)$.) Also, denote by $y^n_{k,j}(t)$ the fraction of queues that are in pool $j$ 
and have queue size {\em exactly} $k \in \barZ_+$:
\beql{eq-y-def}
y^n_{k,j}(t) \doteq x^n_{k,j}(t) - x^n_{k+1,j}(t), ~~k\in \Z_+,
\end{equation}
\beql{eq-y-def-inf}
y^n_{\infty,j}(t) \doteq x^n_{\infty,j}(t) = \lim_{k\to\infty} x^n_{k,j}(t).
\end{equation}

\subsection{Monotonicity.}
\label{subsec-monotone}

The following monotonicity property holds. 
(For a general definition and discussion of stochastic processes' monotonicity 
properties, cf. \cite{Liggett-book}.)
It is a rather straightforward generalization of Lemma 3 in \cite{St2014_pull},
where it was given for the special case of single router. (We give a proof for completeness.)

\begin{lem}
\label{lem-monotone}
Consider two version of the process, $S^n(\cdot)$ and $\bar S^n(\cdot)$ 
[resp. $s^n(\cdot)$ and $\bar s^n(\cdot)$],
with fixed initial states $S^n(0) \le \bar S^n(0)$ [resp. $s^n(0) \le \bar s^n(0)$].
Then, the processes can be constructed on a common probability space,
so that, w.p.1, $S^n(t) \le \bar S^n(t)$ [resp. $s^n(t) \le \bar s^n(t)$]
for all $t\ge 0$. Consequently, 
$S^n(t) \le_{st} \bar S^n(t)$ [resp. $x^n(t) \le_{st} \bar x^n(t)$]
for all $t\ge 0$. This property still holds if the second process $\bar S^n(\cdot)$
[resp. $\bar s^n(\cdot)$] is for a system with 
possibly larger buffer sizes:
$B_j \le \bar B_j \le \infty, ~\forall j\in \cj$. 
\end{lem}

{\em Proof.} It suffices to prove the result for $S^n(\cdot)$ and $\bar S^n(\cdot)$.
We will refer to the systems, corresponding to $S^n(\cdot)$ and $\bar S^n(\cdot)$,
as ``smaller'' and ``larger'', respectively.
Let us consider PULL-2 algorithm first. 
It is clear how to couple 
the service completions in the two systems,
so that any service completion preserves the
$S^n(t) \le \bar S^n(t)$ condition. We make the arrival process
to each router $r$ common for both systems.
It suffices to show that condition $S^n(t) \le \bar S^n(t)$ is preserved 
after an arrival at a given router $r$ at time $t$. Suppose the (joint) system state just before this arrival is such that $S^n \le \bar S^n$. 
If router $r$ does not have pull-messages in either system, i.e. $\sum_i I\{D_i^n=r\} = \sum_i I\{\bar D_i^n=r\} =0$, 
we make a {\em common} random uniform assignment of the arrival to one of the servers in the entire system.
If router $r$ has pull-messages in the smaller system, but not in the larger system,
i.e. $\sum_i I\{D_i^n=r\} > \sum_i I\{\bar D_i^n=r\} =0$, 
we make {\em independent} assignments in the two systems, according to the algorithm.
In the case when router $r$ has pull-messages in both systems,
i.e. $\sum_i I\{D_i^n=r\} \ge \sum_i I\{\bar D_i^n=r\} >0$, 
the corresponding set of idle servers  in the larger system
forms a subset of that in the smaller one.
In this case we do the following.
We make uniform random choice of a pull-message in the smaller system, and 
assign the arrival to the corresponding idle server; if the chosen pull-message in the smaller system is also present
in the larger system, we assign the arrival to the same idle server in the larger system;
if the chosen pull-message in the smaller system is not present
in the larger system, we do an additional step of choosing randomly uniformly 
a pull-message in the larger system, and assigning the arrival accordingly.
Clearly, condition  $S^n(t) \le \bar S^n(t)$ is preserved in each case, and the procedure conforms 
to the PULL-2 algorithm in both systems.

The coupling construction for the PULL-1 algorithm is same, except it is simpler: if an arrival at router $r$ does not find a pull-message there, it is immediately blocked. The condition $S^n(t) \le \bar S^n(t)$ is still preserved
after each arrival.
$\Box$

Note that the state space has a {\em maximum} (``largest'') state,
namely the ``full'' state, with $Q_i=B_j$ for all $j$ and $i\in \cn_j$.
Then, from Lemma~\ref{lem-monotone} we obtain the following

\begin{cor}
\label{cor-monotone}
Consider the process $S^n(\cdot)$ starting from the maximum initial state. Then the process  
$S^n(\cdot)$ is {\em monotone non-increasing}; namely, for any $0\le t_1 \le t_2 \le \infty$,
$S^n(t_2) \le_{st} S^n(t_1)$.
\end{cor}

We will also need an order relation on a subset of servers. We say that {\em state $S'$ is dominated by state $S''$
on a subset $\widehat \cn \subseteq \cn$ of servers}, if $D' \le D''$ 
(understood as \eqn{eq-D-order}) and $Q'_i \le Q''_i, ~i\in \widehat \cn$. Therefore, this relation is such that
if $i\not\in \widehat\cn$ and $Q''_i > 0$ then condition $Q'_i \le Q''_i$ is {\em not} required.

\subsection{Implications of monotonicity: single-router Vs. multi-router case}

Monotonicity is one of the main tools used in this paper, 
as it was in \cite{St2014_pull} for the special case of a single router.
However, some key parts of the analysis in \cite{St2014_pull} {\em cannot} be extended to the multi-router case.
The principal reason is that in the single-router case the state space has the minimum state -- this gives a 
simple characterization of the stationary distribution as the lower invariant measure. In the multi-router case there is {\em no minimum state}, which substantially complicates the analysis. We now discuss these issues in more detail.

{\bf Single router.} In the special case of a single router, $R=1$, studied in \cite{St2014_pull}, the state space has the {\em minimum} (``smallest'') element, namely the state with all servers idle. In that case, for a given system 
and parameter $n$, we considered 
the process, starting from the minimum state $S^n(0)$. Such process, as easily follows from  Lemma~\ref{lem-monotone},
is stochastically non-decreasing in time
$$
S^n(t_1) \le_{st} S^n(t_2) ~~\mbox{[resp. $s^n(t_1) \le_{st} s^n(t_2)$]}, 
~~~0\le t_1 \le t_2 < \infty,
$$
and since the state space is compact, we have
$$
S^n(t) \Rightarrow S^n(\infty) 
~~\mbox{[resp. $s^n(t) \Rightarrow s^n(\infty)$]}, 
~~~ t\to\infty,
$$
where the distribution of $S^n(\infty)$ [resp. $s^n(\infty)$] 
is the {\em lower invariant measure} of process 
$S^n(\cdot)$ [resp. $s^n(\cdot)$]. (The lower invariant measure
is a stationary distribution of the process, stochastically dominated
by any other stationary distribution. Cf. \cite{Liggett-book},
in particular Proposition I.1.8(d).)  Then, as shown in \cite{St2014_pull},
the following is true for
the process $S^n(\cdot)$ [resp. $s^n(\cdot)$], 
{\em as originally defined (without infinite queues)}:
the process is ergodic if and only if $S^n(\infty)$ [resp. $s^n(\infty)$] 
is proper in the sense that 
$$
\pr\{Q_i^n(\infty)<\infty,~\forall i\} = 1~~~ 
\mbox{[resp. $\pr\{x^n_{\infty,j}(\infty)=0,~\forall j\} = 1$]}.
$$
And if the original process is ergodic, the lower invariant measure is its unique
stationary distribution. 

{\bf Multiple routers.} The more general model that we consider in this paper, with multiple routers, $R\ge 1$, is substantially
different. There is a finite number (greater than one) of minimal states; {\em a single minimum state does not exist.} There exists the {\em maximum} (``largest'') state,
namely the ``full'' state, with $Q_i=B_j$ for all $j$ and $i\in \cn_j$. In the special case when
{\em all buffer sizes are finite}, $B_j<\infty, ~\forall j$, the situation is in a sense analogous to that
for the single router in \cite{St2014_pull}. First, since the state space is finite the process is automatically ergodic.
Second, we can consider the process starting with the maximum state; by Corollary~\ref{cor-monotone}
it is stochastically monotone non-increasing, and the distribution of $S^n(t)$ converges to the {\em upper invariant measure}, which is the
unique stationary distribution. However, if buffer sizes are infinite in some or all pools, the situation is much more complicated. The maximum state and the upper invariant measure still exist, but this upper invariant measure is not necessarily (and typically is not) the stationary distribution of the original process. For example, if 
all buffer sizes are infinite, $B_j=\infty$ for all $j$, then the maximum state (with all queues infinite) is invariant, 
so the upper invariant measure is concentrated on it; but it is not the stationary
distribution the original process, if the latter happens to be ergodic. Thus, the upper invariant measure is not useful in the case of infinite buffers, which means that a key element of the approach in \cite{St2014_pull} does not 
extend to multiple routers.

\section{Fluid limits}
\label{sec-fluid-many-server}

In this section we consider limits of the sequence of process
$s^n(\cdot)$ as $n\to\infty$. 
First, define fluid sample paths (FSP), 
which arise as limits of the 
trajectories $s^n(\cdot)$
as $s\to\infty$. The FSP construction will apply to both PULL-1 and PULL-2; it will be clear which parts of the construction become trivial (or redundant) for either PULL-1 or PULL-2. Similarly, when we establish properties of the FSPs,
they hold for both PULL-1 and PULL-2, unless explicitly stated otherwise.

Without loss of generality, assume that
 Markov process $s^n(\cdot)$ for each $n$ is driven by a
common set of primitive processes, as defined next.

Let $A_r^n(t), ~t\ge 0$, denote the number of exogenous arrivals into the system, via router $r$,
in the interval $[0,t]$. Assume that 
\beql{eq-driving-arr}
A_r^n(t) = \Pi_r^{(a)}(\frac{\lambda}{R} n t),
\end{equation}
where $\Pi_r^{(a)}(\cdot)$ is an independent unit rate Poisson process.
The functional strong law of large numbers (FSLLN) holds:
\beql{eq-flln-poisson-arr}
\frac{1}{n}\Pi_r^{(a)}(nt) \to t, ~u.o.c., ~~w.p.1, ~~\forall r\in\calr.
\end{equation}
Denote by $D^n_{k,j,r}(t), ~t\ge 0, ~1 \le k < \infty$, 
the total number of departures in $[0,t]$
from servers in pool $j$ with queue length $k$, which are ``associated'' with router $r$;
this association with router $r$ is only used when $k=1$, in which case an $r$-pull-message is generated.
(Also note that in case of PULL-1, $D^n_{k,j,r}(t)\equiv 0$ for any $k\ge 2$.)
Assume that
\beql{eq-driving-dep}
D^n_{k,j,r}(t) = \Pi^{(d)}_{k,j,r} \left(\frac{1}{R}\int_0^t n y_{k,j}^n(w) \mu_j dw\right),
\end{equation}
where $\Pi^{(d)}_{k,j,r}(\cdot)$ are independent unit rate Poisson processes.
(Recall that departures from -- and arrivals to -- infinite queues can be ignored,
in the sense that they do not change the system state.)
Similarly to \eqn{eq-flln-poisson-arr}, 
\beql{eq-flln-poisson-dep}
\frac{1}{n}\Pi^{(d)}_{k,j,r}(nt) \to t, ~u.o.c., ~~w.p.1, ~~~ 1 \le k < \infty, ~\forall j,~\forall r\in\calr.
\end{equation}

The routing of new arrivals, going through router $r$, is constructed as follows.
There are two sequences of i.i.d. random variables, 
$$
\chi_r(1), \chi_r(2), \ldots, ~~~~\mbox{and}~~~ \zeta_r(1), \zeta_r(2), \ldots,
$$
uniformly distributed in $[0,1)$. The routing of the $m$-th arrival (via router $r$) into the system
is determined by the values of r.v. $\chi_r(m)$ and $\zeta_r(m)$,
as follows. (We will drop index $m$, because we consider one arrival.)
Let $s^n$ denote the system state just before the arrival.
If $\sum_j \xi^n_{r,j} = 0$, i.e. there are no pull-messages at $r$, 
then under PULL-1 the arrival is blocked, and under PULL-2
the routing is determined
by $\zeta_r$ as follows.
The customer is sent to a server with $k$, $k\ge 0$, customers in pool 1,
if $\zeta_r \in [x^n_{k+1,1},x^n_{k,1})$, and to a server with $k=\infty$ customers in pool 1, if
$\zeta_r \in [0,x^n_{\infty,1})$; the customer is sent to a server with $k$, $k\ge 0$, customers in pool 2,
if $\zeta_r \in [\beta_1+x^n_{k+1,2},\beta_1 + x^n_{k,2})$, 
and to a server with $k=\infty$ customers in pool 2, if
$\zeta_r \in [\beta_1, \beta_1 + x^n_{\infty,2})$; and so on.
If $\sum_j \xi^n_{r,j} > 0$, i.e. there are pull-messages at $r$,  the routing is determined
by $\chi_r$ as follows.
Let $a=\sum_j \xi^n_{r,j}$, $p_j = \xi^n_{r,j}/a$.
If $\chi_r \in [0,p_1)$, the customer is routed to an idle server associated with $r$ in pool 1; 
if $\chi_r \in [p_1,p_1+p_2)$ -- it is routed to an idle server associated with $r$ in pool 2; and so on.

Denote
$$
f^n_r(w,u) \doteq \frac{1}{n} \sum_{m=1}^{\lfloor nw \rfloor} I\{\chi_r(m) \le u\}, ~~
g^n_r(w,u) \doteq \frac{1}{n} \sum_{m=1}^{\lfloor nw \rfloor} I\{\zeta_r(m) \le u\}, ~~
$$
where $w\ge 0$, $0\le u < 1$.  Obviously, from the strong law of large numbers
 and the monotonicity of $f^n(w,u)$ and $g^n(w,u)$ on both arguments, we 
have the following FSLLN:
\beql{eq-flln-random}
f^n_r(w,u) \to wu, ~~g^n_r(w,u) \to wu, 
~~~\mbox{u.o.c.}, ~~~~~w.p.1, ~~\forall r\in \calr.
\end{equation}
It is easy (and standard) to see that, for any $n$, w.p.1,
 the realization of the process $s^n(\cdot)$ 
is uniquely determined by the
initial state $s^n(0)$
and the realizations of the driving processes
$\Pi_r^{(a)}(\cdot)$, $\Pi^{(d)}_{k,j,r}(\cdot)$, $\chi_r(\cdot)$ and $\zeta_r(\cdot)$.

Denote by $A^n_{k,j,r}(t), ~k\in \Z_+, ~t\ge 0$, 
the total number of arrivals  in $[0,t]$ via router $r$
into servers in pool $j$ with queue length $k$. 
(Recall that arrivals to infinite queues can be ignored.
Also note that, under PULL-1, $A^n_{k,j,r}(t)\equiv 0$ if $k\ge 1$.)
Obviously, for any $0 \le t_1 \le t_2 < \infty$
\beql{eq-200star}
\sum_j \sum_{0 \le k < \infty} [A^n_{k,j,r}(t_2)-A^n_{k,j,r}(t_1)] \le A^n_r(t_2)-A^n_r(t_1),
~~~r\in\calr.
\end{equation}
(Note that, under PULL-1,
$$
[A^n_r(t_2)-A^n_r(t_1)] - \sum_j [A^n_{0,j,r}(t_2)-A^n_{0,j,r}(t_1)],
$$
is the number of arrivals via $r$ in $(t_1,t_2]$ which are blocked.)

We define the fluid-scaled the arrival and departure processes:
$$
a^n_{k,j,r}(t) = \frac{1}{n} A^n_{k,j,r}(t),~~ 0 \le k < \infty,
$$
$$
d^n_{k,j,r}(t) = \frac{1}{n} D^n_{k,j,r}(t),~~~ 1 \le k < \infty.
$$

\begin{definition}[Fluid sample path]
\label{def-fsp}
A set of uniformly Lipschitz continuous functions
$s(\cdot)=[x_{k,j}(\cdot),~ k\in \Z_+, ~j\in \cj; ~ \xi_{r,j}(\cdot),~ r\in \calr, ~j\in \cj]$
on the time interval $[0,\infty)$ we call a {\em fluid sample path} (FSP), if there exist
realizations of the primitive driving processes,
 satisfying conditions \eqn{eq-flln-poisson-arr}, \eqn{eq-flln-poisson-dep}
and \eqn{eq-flln-random}
and a fixed subsequence of $n$, along which
\beql{eq-fsp-def}
s^n(\cdot) \to s(\cdot), ~~~u.o.c.,
\end{equation}
and also
\beql{eq-fsp-def2}
a^n_{k,j,r}(\cdot) \to a_{k,j,r}(\cdot), ~~u.o.c., ~~~0\le k < \infty, ~j\in\cj, ~r\in\calr,
\end{equation}
\beql{eq-fsp-def3}
d^n_{k,j,r}(\cdot) \to d_{k,j,r}(\cdot), ~~u.o.c., ~~~1\le k < \infty, ~j\in\cj, ~r\in\calr,
\end{equation}
for some functions $a_{k,j,r}(\cdot)$ and $d_{k,j,r}(\cdot)$.
\end{definition}

Given the metric \eqn{eq-metric}
on $\cs$, condition \eqn{eq-fsp-def} is equivalent to
component-wise convergence:
$$
x^n_{k,j}(\cdot) \to x_{k,j}(\cdot), ~~~u.o.c., ~~k\in \Z_+, ~j\in \cj,
$$
$$
\xi^n_{r,j}(\cdot) \to \xi_{r,j}(\cdot), ~~~u.o.c., ~~r\in \calr, ~j\in \cj.
$$
Also, since all functions $a_{k,j,r}^n(\cdot)$ and $d_{k,j,r}^n(\cdot)$ are non-decreasing, the FSP definition
easily implies that all limit functions $a_{k,j,r}(\cdot)$ and $d_{k,j,r}(\cdot)$ are necessarily non-decreasing, uniformly Lipschitz continuous, and such that $a_{k,j,r}^n(0)=0$, $d_{k,j,r}^n(0)=0$.

In what follows, WLOG, we use the convention that each FSP $s(\cdot)$ has an associated set 
of functions $a_{k,j,r}(\cdot)$ and $d_{k,j,r}(\cdot)$, which are a part of its definition.
(A cleaner, but more cumbersome FSP definition should include functions 
$a_{k,j,r}(\cdot)$ and $d_{k,j,r}(\cdot)$ as FSP components. We do not do this to simplify notation.)

For any FSP $s(\cdot)$, almost all points $t\ge 0$ (w.r.t. Lebesgue measure)
are {\em regular}; namely, all component functions 
(including $a_{k,j,r}(\cdot)$ and $d_{k,j,r}(\cdot)$)
have proper 
(equal right and left) derivatives.
Note that $t=0$ is {\em not} a regular point.

Analogously to notation in \eqn{eq-x-inf} - \eqn{eq-y-def-inf}, we will denote:
$$
x_{\infty,j}(t)\doteq \lim_{k\to\infty} x_{k,j}(t)
$$
$$
y_{k,j}(t) \doteq x_{k,j}(t) - x_{k+1,j}(t), ~~k\in \Z_+,
$$
$$
y_{\infty,j}(t) \doteq x_{\infty,j}(t) = \lim_{k\to\infty} x_{k,j}(t).
$$
For two FSPs $s(\cdot)$ and $\bar s(\cdot)$,
$s(\cdot) \le \bar s(\cdot)$ is defined as $s(t) \le \bar s(t), ~t\ge 0$.

\begin{lem}
\label{lem-ms-fsp-conv}
Consider a sequence in $n$ of processes $s^n(\cdot)$ with deterministic initial states
$s^n(0) \to s(0)\in \cs$.
Then w.p.1 any subsequence of $n$ has a further subsequence, along which
$$
s^n(\cdot) \to s(\cdot) ~~~u.o.c.,
$$
where $s(\cdot)$ is an FSP.
\end{lem}

{\em Proof.} The proof  is analogous to that of Lemma 5 in \cite{St2014_pull}.
We give it here for completeness.
All processes $a^n_{k,j,r}(\cdot)$ and $d^n_{k,j,r}(\cdot)$ are non-decreasing.
W.p.1 the primitive processes satisfy the FSLLN \eqn{eq-flln-poisson-arr},
\eqn{eq-flln-poisson-dep} and \eqn{eq-flln-random}. 
From here it is easy to observe the following:
w.p.1 any subsequence of $n$ has a further subsequence along which
the u.o.c. convergences
$$
a^n_{k,j,r}(\cdot) \to a_{k,j,r}(\cdot), ~~d^n_{k,j,r}(\cdot) \to d_{k,j,r}(\cdot),
$$
hold for all $(k,j,r)$, with the limiting functions
$a_{k,j,r}(\cdot)$ and $d_{k,j,r}(\cdot)$ being non-decreasing,
uniformly Lipschitz continuous. 
Then, we use obvious relations between pre-limit process components:
\beql{eq-prel-basic1}
x^n_{1,j}(t) - x^n_{1,j}(0) = \sum_r a^n_{0,j,r}(t) - \sum_r d^n_{1,j,r}(t) - \sum_r a^n_{1,j,r}(t),
~~~j\in\cj,
\end{equation}
\beql{eq-prel-basic2}
x^n_{k,j}(t) - x^n_{k,j}(0) = \sum_r a^n_{k-1,j,r}(t)  - \sum_r d^n_{k,j,r}(t), ~~~j\in\cj, ~2\le k <\infty,
\end{equation}
\beql{eq-prel-basic3}
\xi^n_{r,j}(t) - \xi^n_{r,j}(0) = d^n_{1,j,r}(t) - a^n_{0,j,r}(t), ~~~j\in\cj, ~ r\in\calr.
\end{equation}
The result easily follows. 
$\Box$

In the following Lemmas \ref{lem-ms-fsp-prop-basic} and \ref{lem-ms-fsp-prop}
we establish some properties of the FSPs; first, very basic (in Lemma \ref{lem-ms-fsp-prop-basic})
and then more involved (in Lemma \ref{lem-ms-fsp-prop}).
These properties hold for both PULL-1 and PULL-2.
(As stated at the beginning of this section, unless specified otherwise, all FSP properties that we prove,
hold for both algorithm versions.) This does {\em not} mean, of course, that FSP behavior is same under 
PULL-1 and PULL-2. The properties in Lemmas \ref{lem-ms-fsp-prop-basic} and \ref{lem-ms-fsp-prop}
 do {\em not} (and not intended to)
characterize FSP dynamics completely; but they will suffice for the proofs of the main results of this paper.
Also note that, under PULL-1, for any FSP, $x_{k,j}(t) = y_{k,j}(t) \equiv 0, ~t\ge 0,$ for any $k \ge 2$. This means, for example, that, for PULL-1, \eqn{eq-deriv-2} becomes trivial and property (vii) in
Lemma -- irrelevant (because the 'if' condition never holds).

\begin{lem}
\label{lem-ms-fsp-prop-basic}
Any FSP $s(\cdot)$ has the following properties.\\
(i) For any $t\ge 0$,
\beql{eq-fl-basic1}
x_{1,j}(t) - x_{1,j}(0) = \sum_r a_{0,j,r}(t) - \sum_r d_{1,j,r}(t) - \sum_r a_{1,j,r}(t),
~~~j\in\cj,
\end{equation}
\beql{eq-fl-basic2}
x_{k,j}(t) - x_{k,j}(0) = \sum_r a_{k-1,j,r}(t)  - \sum_r d_{k,j,r}(t), ~~~j\in\cj, ~2\le k <\infty,
\end{equation}
\beql{eq-fl-basic3}
\xi_{r,j}(t) - \xi_{r,j}(0) = d_{1,j,r}(t) - a_{0,j,r}(t), ~~~j\in\cj, ~ r\in\calr,
\end{equation}
\beql{eq-fl-basic4}
d_{k,j,r}(t) = \frac{1}{R} \int_0^t \mu_j y_{k,j}(u)du, ~~~2\le k <\infty, ~ j\in\cj, ~ r\in\calr.
\end{equation}
(ii) For any $0 \le t_1 \le t_2$, and any $r\in \calr$,
\beql{eq-fl-basic5}
\sum_{k,j} a_{k,j,r}(t_2) - \sum_{k,j} a_{k,j,r}(t_1) \le \frac{\lambda}{R} (t_2 - t_1).
\end{equation}
(iii) If for some $r\in \calr$, 
\beql{eq-idleness-basic}
\sum_{j} \xi_{r,j}(t) > 0, ~~ t \in [t_1,t_2],
\end{equation}
then:
\beql{eq-fl-basic6}
a_{k,j,r}(t_2) - a_{k,j,r}(t_1) = 0, ~~k\ge 1, ~j \in \cj,
\end{equation}
\beql{eq-fl-basic7}
a_{0,j,r}(t_2) - a_{0,j,r}(t_1) = 
\frac{\lambda}{R} \int_{t_1}^{t_2} \frac{\xi_{r,j}(t)}{\sum_\ell \xi_{r,\ell}(t)} dt, ~~j \in \cj,
\end{equation}
\beql{eq-fl-basic8}
\sum_j a_{0,j,r}(t_2) - \sum_j a_{0,j,r}(t_1) = \frac{\lambda}{R} (t_2 - t_1).
\end{equation}
\end{lem}

{\em Proof.} (i) Properties \eqn{eq-fl-basic1}-\eqn{eq-fl-basic4} follow from the corresponding relations 
for the pre-limit process trajectories, \eqn{eq-prel-basic1}-\eqn{eq-prel-basic3}, \eqn{eq-driving-dep}, and the FSP definition.

(ii) Similarly to (i), \eqn{eq-fl-basic5} follows from \eqn{eq-200star}, \eqn{eq-driving-arr}, and the FSP definition.

(iii) If condition \eqn{eq-idleness-basic} holds, 
then a sequence $s^n(\cdot)$ of pre-limit trajectories defining FSP $s(\cdot)$,
is such that for all sufficiently large $n$ the pre-limit version of \eqn{eq-idleness-basic} holds as well:
$$
\sum_{j} \xi^n_{r,j}(t) > 0, ~~ t \in [t_1,t_2].
$$
This means that (for all large $n$) the trajectory $s^n(\cdot)$ 
is such that in the interval $[t_1,t_2]$ all new arrivals
 are routed to idle servers;
furthermore, in a small neighborhood of time $t\in [t_1,t_2]$ 
the fraction of arrivals to router $r$ that are sent to pool $j$ is close to $\xi_{r,j}(t)/[\sum_{\ell} \xi_{r,\ell}(t)]$.
This and FSP definition easily imply \eqn{eq-fl-basic6} and \eqn{eq-fl-basic7} -- we omit formalities,
which are rather standard. Property \eqn{eq-fl-basic8} follows from \eqn{eq-fl-basic7}.
$\Box$

\begin{lem}
\label{lem-ms-fsp-prop}
(i) Suppose two FSPs, $s(\cdot)$ and $\bar s(\cdot)$, are such that $s(0) \le \bar s(0)$,
and for some $\tau>0$ (including possibly $\tau=\infty$)
 these FSPs are the unique in $[0,\tau)$ for the initial states 
$s(0)$ and  $\bar s(0)$, respectively. Then $s(t) \le \bar s(t)$ for $0\le t < \tau$.
\\ (ii) For any FSP $s(\cdot)$ the following holds. If $t\ge 0$ is such that 
\beql{eq-idleness}
\sum_j \xi_{r,j}(t) >0, \forall r\in\calr,
\end{equation}
then the derivatives
$(d/dt) x_{k,j}(t)$ and $(d/dt) \xi_{r,j}(t)$ exist (if this $t=0$ -- right derivatives exist)
and equal to
\beql{eq-deriv-1}
\frac{d}{dt}  x_{1,j}(t) =  \frac{\lambda}{R}\sum_r \frac{ \xi_{r,j}(t)}{\sum_{\ell}  \xi_{r,\ell}(t)} - \mu_j  y_{1,j}(t),~ j\in \cj,
\end{equation}
\beql{eq-deriv-2}
\frac{d}{dt}  x_{k,j}(t) = - \mu_j  y_{k,j}(t) \le 0, ~~ 2\le k < \infty, ~ j\in \cj,
\end{equation}
\beql{eq-deriv-3}
\frac{d}{dt} \xi_{r,j}(t) = 
\frac{1}{R}  y_{1,j} \mu_j(t) - \frac{\lambda}{R} \frac{\xi_{r,j}(t)}{\sum_{\ell} \xi_{r,\ell}(t)}, ~r\in\calr, ~j\in \cj.
\end{equation}
\\ (iii) If initial state $s(0)$ of an FSP is such that \eqn{eq-idleness} holds (with $t=0$)
and $\sum_j x_{2,j}(0)=0$,
then the FSP is unique in the interval $[0,\tau)$, where $\tau$ is the smallest time $t$
when \eqn{eq-idleness} no longer holds; $\tau=\infty$ if such $t$ does not exist.
\\ (iv) An FSP $s(\cdot)$ with initial condition $s(0)=s^*$ is unique, 
and it is stationary, i.e. $s(t) \equiv s^*$.
\\ (v) Suppose, for an FSP $s(\cdot)$ and some $t\ge 0$ we have:
$\xi_{r,j}(t)=0,~\forall r, ~ \forall j$; $y_{1,j}(t) \mu_j>0, ~\forall j$;
 $\sum_j y_{1,j}(t) \mu_j>\lambda$. Then, the following right derivatives exist and equal to
\beql{eq-deriv-init}
\frac{d^+}{dt} \xi_{r,j}(t) = \frac{1}{R} \left[ y_{1,j}(t) \mu_j - \frac{y_{1,j}(t) \mu_j}{\sum_{\ell} y_{1,\ell}(t) \mu_{\ell}}\lambda  \right] > 0, ~r\in\calr, ~j\in \cj.
\end{equation}
\\ (vi) Consider the system with $B_j = 1$ for all $j$, and an FSP $\bar s(\cdot)$ 
for this system,
starting from the maximum initial state, that is
$\bar x_{1,j}(0) = \beta_j, ~ \bar x_{2,j}(0) = 0,$ for all $j$. Such an FSP $\bar s(\cdot)$ is unique, 
 monotone non-increasing, $\bar s(t_1)\ge \bar s(t_2), ~t_1 \le t_2$,
and is such that $\bar s(t) \to s^*$ as $t\to\infty$. 
\\ (vii) For any $\epsilon>0$, there exist $\tau>0$ and $\delta>0$,
such that the following holds for any FSP $s(\cdot)$.
If at time $t\ge 0$, $s(t) \ge s^*$, 
$x_{1,j}(t)=\nu_{j}$ for all $j\in \cj$,
and $x_{2,\ell}(t)\ge \epsilon$ for some fixed $\ell$,
then
$$
x_{1,\ell}(\tau) \ge \nu_{\ell}+\delta.
$$
\end{lem}

\begin{rem}
\label{rem-deriv} 
Note that, although the Lipschitz continuity of an FSP components guarantees the existence of their derivatives
almost everywhere, as opposed to everywhere, statement (ii) of Lemma~\ref{lem-ms-fsp-prop} says that the derivatives 
\eqn{eq-deriv-1}-\eqn{eq-deriv-3} do in fact exist at a time $t$ where \eqn{eq-idleness} holds. Similarly,
the right derivatives \eqn{eq-deriv-init} do exist at a time $t$ at which the conditions specified in statement (v) of Lemma~\ref{lem-ms-fsp-prop} hold.
\end{rem}

{\em Proof of Lemma~\ref{lem-ms-fsp-prop}.} (i) Consider the sequences of {\em processes} $s^n(\cdot)$ and $\bar s^n(\cdot)$,
with fixed initial states such that $s^n(0) \le \bar s^n(0)$, and $s^n(0)\to s(0), ~\bar s^n(0)\to \bar s(0)$
as $n\to\infty$. By Lemma~\ref{lem-monotone}, $s^n(t) \le_{st} \bar s^n(t)$ for all $t\ge 0$.
This, Lemma~\ref{lem-ms-fsp-conv}, and the FSPs uniqueness in $[0,\tau)$ imply the result.

(ii) If \eqn{eq-idleness} holds for $t\ge 0$, then \eqn{eq-idleness-basic} holds in a small neighborhood 
of $t$ for all $r\in \calr$. Then we apply properties (iii) and (i) of Lemma~\ref{lem-ms-fsp-prop-basic}
to easily obtain the existence of -- and the expressions for -- the derivatives 
\eqn{eq-deriv-1}-\eqn{eq-deriv-3} at $t$.

(iii) From (ii) we in particular have the following. For an FSP $s(\cdot)$,
at any $t\ge 0$ such that conditions \eqn{eq-idleness} and 
\beql{eq-idleness2}
x_{2,j}(t)=0~~\mbox{[or, equivalently, $y_{1,j}(t)=x_{1,j}(t)$]},~~~\forall j,
\end{equation}
hold, we have
$$
(d/dt) x_{k,j}(t) = 0, ~~k\ge 2, ~\forall j.
$$
Therefore, condition \eqn{eq-idleness2} must hold in the interval $[0,\tau)$.
In that interval, $x_{1,j}(t) = y_{1,j}(t) = \beta_j - \sum_r \xi_{r,j}(t)$, and therefore \eqn{eq-deriv-3} becomes
\beql{eq-deriv-33}
\frac{d}{dt} \xi_{r,j}(t) = 
\frac{1}{R}  (\beta_j - \sum_{r'} \xi_{r',j}(t))\mu_j - \frac{\lambda}{R} \frac{\xi_{r,j}(t)}{\sum_{\ell} \xi_{r,\ell}(t)}, ~r\in\calr, ~j\in \cj.
\end{equation}
We see that, in $[0,\tau)$, the FSP trajectory is uniquely determined by the solution to ODE \eqn{eq-deriv-33}.

(iv) By (ii) and the definition of $x^*$, $(d/dt)s(t)=0$ if $s(t)=s^*$.
Then we apply (iii).

(v) The proof will be in two steps. Since the statement is about right derivatives at $t$, WLOG, to simplify notation,
 let $t=0$ and we consider an FSP with time being $u\ge 0$.

Step 1. We will use notation $w_j(u)=y_{1,j}(u) \mu_j$.
In this step we prove that, for any $r$, in a sufficiently small neighborhood $[0,\epsilon]$,
\beql{eq-right-deriv}
\frac{d}{du} \sum_j \xi_{r,j}(u) = \frac{1}{R} 
\sum_{j} w_j(u)   - \frac{1}{R} \lambda, ~~u\in (0,\epsilon],
\end{equation}
\beql{eq-right-deriv-0}
\frac{d^+}{du} \sum_j \xi_{r,j}(0) = \frac{1}{R} 
\sum_{j} w_j(0)   - \frac{1}{R} \lambda.
\end{equation}
The meaning of \eqn{eq-right-deriv}-\eqn{eq-right-deriv-0} is simple: the first term in the RHS is the (fluid-scaled) rate at which new $r$-pull-messages are generated, and the second term is the rate at which they are ``used'' (and removed). Note that since the second term is equal to $\lambda/R$ for all $r$, \eqn{eq-right-deriv}-\eqn{eq-right-deriv-0} mean that in $[0,\epsilon]$ ``all'' (in the fluid limit) new arrivals go to idle servers.
Note also that it suffices to prove \eqn{eq-right-deriv}, because \eqn{eq-right-deriv-0} follows from  
\eqn{eq-right-deriv} by continuity of $w_j(\cdot)$ (and the fact that $\sum_j \xi_{r,j}(\cdot)$
is absolutely continuous).

The formal argument is as follows. Choose $\epsilon>0$ small enough so that the RHS of
\eqn{eq-right-deriv} is positive in $[0,\epsilon]$. (It is positive for $u=0$, and $w_j(\cdot)$ are continuous.) 
Almost everywhere in $[0,\epsilon]$,
$$
\frac{d}{du} \sum_{r,j} \xi_{r,j}(u) = 
\sum_{r,j} \frac{d}{du} d_{1,j,r}(u) - \sum_{r,j} \frac{d}{du} a_{0,j,r}(u) \ge 
$$
$$
\sum_{r,j} \frac{1}{R} w_j(u) - \lambda = \sum_{j} w_j(u) - \lambda >0,
$$
where we used Lemma~\ref{lem-ms-fsp-prop-basic}(i)-(ii). Threrefore, $\sum_{r,j} \xi_{r,j}(u)>0$ in $(0,\epsilon]$.
Pick any $u'\in (0,\epsilon]$. There
exists $r'$ for which $\sum_j \xi_{r',j}(u) > 0$ in a small interval $(u'-\epsilon',u'+\epsilon')$. 
Consider now the derivative (it exists a.e. in $(u'-\epsilon',u'+\epsilon')$):
$$
\frac{d}{du} \sum_{r\ne r'} \sum_{j} \xi_{r,j}(u) = 
\sum_{r\ne r'} \sum_{j} \frac{d}{du} d_{1,j,r}(u) - \sum_{r\ne r'} \sum_{j} \frac{d}{du} a_{0,j,r}(u) \ge 
$$
$$
\sum_{r\ne r'} \sum_{j} \frac{1}{R} w_j(u) - \sum_{r\ne r'} \sum_{j} \frac{1}{R} \lambda =
\frac{R-1}{R} [\sum_{j} w_j(u) - \lambda]>0.
$$
We conclude that $\sum_{r\ne r'} \sum_{j} \xi_{r,j}(u)>0$ in $(u'-\epsilon',u'+\epsilon')$,
and therefore there exists $r''$ such that $\sum_j \xi_{r'',j}(u) > 0$ in a small interval $(u'-\epsilon'',u'+\epsilon'')$.
Continuing this argument, we obtain that
$\sum_j \xi_{r,j}(u) > 0, \forall r,$ for $u=u'$, and therefore for all $u \in (0,\epsilon]$. In other words,
\eqn{eq-idleness} holds in $(0,\epsilon]$, so that \eqn{eq-deriv-3} holds as well. Summing up \eqn{eq-deriv-3} over
$j$, we obtain \eqn{eq-right-deriv}. This proves Step 1.

Step 2. Now we will establish \eqn{eq-deriv-init}, which, using the notations  of this proof, can be written as
\beql{eq-y0j}
\frac{d^+}{du} \xi_{r,j}(0) = b_j \equiv \frac{1}{R}
\left[ w_j(0) - \frac{w_j(0)}{\sum_{\ell} w_\ell(0)}\lambda  \right]
= \frac{w_j(0)}{R}
\left[ \frac{\sum_{\ell} w_\ell(0) -\lambda}{\sum_{\ell} w_\ell(0)}  \right].
\end{equation}
It follows from \eqn{eq-right-deriv} (and it is already established in the proof of Step 1), that
in a sufficiently small interval $(0,\epsilon]$ 
 condition \eqn{eq-idleness} holds,  and then \eqn{eq-deriv-3} holds as well. In addition, in 
$(0,\epsilon]$, $w_j(u) >0$ for any $j$. Then, from \eqn{eq-deriv-3} we see that $\xi_{r,j}(u)>0$ for
$u\in (0,\epsilon]$, any $j$ and any $r$. (Indeed, if $\xi_{r,j}(u)=0$ for some $u\in (0,\epsilon]$, then by 
\eqn{eq-deriv-3}
$(d/du)\xi_{r,j}(u)>0$ at that $u$, which is impossible.)

Fix any $j$ and $r$. 
Consider function 
$$
g_{r,j}(u) = \ln [\xi_{r,j}(u) / (b_j u)].
$$
For $u\in (0,\epsilon]$, 
$$ 
\frac{d}{du} \xi_{r,j}(u) =  \frac{1}{R} w_j(u)  -
 \frac{\lambda}{R}  \frac{\xi_{r,j}(u)}{\sum_{\ell} \xi_{r,\ell}(u)}, ~j\in \cj.
$$ 
We can write
$$
g'_{r,j}(u) = \frac{\xi'_{r,j}(u)}{\xi_{r,j}(u)} - \frac{1}{u} = \frac{\xi'_{r,j}(u) - \xi_{r,j}(u)/u}{\xi_{r,j}(u)} =
$$
\beql{eq-g-deriv}
\frac{\frac{1}{R} w_j(u)  -\frac{\lambda}{R}  \frac{\xi_{r,j}(u)/u}{\sum_{\ell} \xi_{r,\ell}(u)/u} - \xi_{r,j}(u)/u}{\xi_{r,j}(u)}.
\end{equation}
In the rest of Step 2 proof, we will write $o(1)$ to signify any function $\gamma(u), ~u\ge 0,$ such that 
$\gamma(u) \to \gamma(0) = 0$ as $u\downarrow 0$; same expression $o(1)$ may signify different such 
functions, even within same formula. We have $w_j(u) = w_j(0) + o(1)$ and, by \eqn{eq-right-deriv-0},
$$
R \sum_{\ell} \xi_{r,\ell}(u)/u = \sum_\ell w_\ell(0) - \lambda + o(1).
$$
Using these formulas in \eqn{eq-g-deriv}, along with the fact that 
$\xi_{r,j}(u)/u$ is bounded (because $\xi_{r,j}(\cdot)$ is Lipschitz), we can rewrite
\beql{eq-g-deriv2}
g'_{r,j}(u)  = 
\frac{\frac{1}{R}w_j(0)[\sum_{\ell} w_\ell(0) - \lambda]   - [\sum_{\ell} w_\ell(0) ]\xi_{r,j}(u)/u + o(1)}    {[\sum_\ell w_\ell(0) - \lambda + o(1)]\xi_{r,j}(u)}.
\end{equation}
From the definition of $b_j$ in \eqn{eq-y0j},
$$
\frac{1}{R}w_j(0)[\sum_{\ell} w_\ell(0) - \lambda] = b_j \sum_{\ell} w_\ell(0).
$$
Substituting this into the first term of the numerator of the RHS of \eqn{eq-g-deriv2}, we obtain
\beql{eq-asypm}
g'_{r,j}(u)  = 
\frac{[1-\frac{\xi_{r,j}(u)}{b_j u}] b_j \sum_\ell w_\ell(0) + o(1)}
{[\sum_\ell w_\ell(0) - \lambda + o(1)]\xi_{r,j}(u)}.
\end{equation}
Since FSP is Lipschitz, $\xi_{r,j}(u) \le c u, ~u\ge 0$, for some $c>0$. Consequently, $g_{r,j}(u) \le c'$ and 
$\xi_{r,j}(u) / (b_j u) \le c'$ for all $u>0$. 
We obtain the following property. For any $\delta>0$, we can choose
$\epsilon>0$ small enough, so that for any $u\in (0,\epsilon]$, 
\beql{eq-cond-deriv}
g_{r,j}(u) \ge \delta ~ (\mbox{i.e.}~ \xi_{r,j}(u)/[b_j u] \ge e^\delta > 1) ~~\mbox{implies}~~ 
g'_{r,j}(u) \le -c_1 /(c_2 u), 
\end{equation}
for some positive constants $c_1$ and $c_2$.
We claim that property \eqn{eq-cond-deriv} implies
\beql{eq-claim-delta}
g_{r,j}(u) \le \delta, ~~ \forall u\in (0,\epsilon].
\end{equation}
Indeed, suppose not, i.e. $g_{r,j}(u) > \delta$
for some $u\in (0,\epsilon]$. Given \eqn{eq-cond-deriv}, it is {\em not} possible
that $g_{r,j}(u') \le \delta$ for some $u' \in (0,u)$, because then $g_{r,j}(\cdot)$ could not exceed $\delta$
in $[u',\epsilon]$.
Therefore, the only remaining possibility is that $g_{r,j}(u') > \delta$
for all $u' \in (0,u)$; but if so, from \eqn{eq-cond-deriv} and $\int_0^u (1/u') du' = \infty$,
we would obtain $g_{r,j}(u') \uparrow \infty$ as $u'\downarrow 0$, which contradicts to the fact that
$g_{r,j}(\cdot)$ is upper bounded. This proves \eqn{eq-claim-delta}.
Since \eqn{eq-claim-delta} holds for an arbitrarily small $\delta>0$ (with corresponding $\epsilon$), we obtain 
\beql{eq-upper55}
\limsup_{u\downarrow 0} g_{r,j}(u) \le 0, ~~\mbox{or, equivalently,} ~~ \limsup_{u\downarrow 0} \frac{\xi_{r,j}(u)}{ u}
\le b_j.
\end{equation}
However, by \eqn{eq-right-deriv-0},
$$
\lim_{u\downarrow 0} \frac{\sum_j \xi_{r,j}(u)}{ u} = \frac{1}{R} 
\sum_{j} w_j(0)   - \frac{1}{R} \lambda
= \sum_j b_j.
$$
This implies that \eqn{eq-y0j} must hold, which completes  Step 2.

(vi) Consider {\em any} fixed FSP $\bar s(\cdot)$ starting from the specified maximum initial state $\bar s(0)$.
We will first prove the properties described in the statement. (After that we will show uniqueness.)
Conditions $\bar x_{2,j}(t) \equiv 0, ~t\ge 0,$ and $\bar x_{1,j}(t) \equiv y_{1,j}(t), ~t\ge 0,$
hold automatically (because it is an FSP for the system with all $B_j=1$). 
The initial state is such that $\bar y_{1,j}(0) = \beta_j$ for all $j$, so that 
$\sum_j \mu_j \bar y_{1,j}(0) > \lambda$, and therefore the conditions of (v) hold for $t=0$.
From (v) we see that condition \eqn{eq-idleness} holds for $\bar s(\cdot)$ at least in some small interval $(0,\epsilon]$.
Therefore, in $(0,\epsilon)$, the FSP satisfies ODE \eqn{eq-deriv-33}. Recall that, by (iii), \eqn{eq-deriv-33}
determines the FSP uniquely up to a time (possibly infinity) when \eqn{eq-idleness} no longer hold.
Also form (v) we observe that for any sufficiently small $\Delta>0$, and any sufficiently small $t_0>0$
(depending on $\Delta$),
we have $\bar s(t_0 + \Delta) \le \bar s(t_0)$. FSPs $\bar s(t_0 + \Delta + \cdot)$ and  $\bar s(t_0+\cdot)$
are the unique FSPs in the interval $[0,\epsilon - t_0 -\Delta]$, starting from 
$\bar s(t_0 + \Delta)$ and $\bar s(t_0)$, respectively. But then, by (i), 
$\bar s(t + \Delta) \le \bar s(t)$ for all $t\in [t_0, \epsilon-\Delta]$. Given that 
$\Delta$ and $t_0$ can be arbitrarily small, we conclude that $\bar s(\cdot)$ is monotone non-increasing 
in $[0,\epsilon]$ (which, recall, implies that each $\bar \xi_{r,j}(\cdot)$ in non-decreasing).
But then condition \eqn{eq-idleness} holds at time $t=\epsilon$, and therefore 
$\bar s(\cdot)$ must continue to be monotone non-increasing beyond time $\epsilon$, up to infinity,
because \eqn{eq-deriv-33} can never be violated.
We proved that $\bar s(\cdot)$ is monotone non-increasing in $[0,\infty)$.

Since $\bar s(\cdot)$ is non-increasing, as $t\to\infty$, $\bar s(t) \to s^{**}$ for some $s^{**}$.
Let us show that $s^{**} = s^*$. For any sufficiently small $t_0>0$, 
we have $\bar s(t_0) \ge s^*$, and the FSPs $\bar s(t_0+\cdot)$ and $s(\cdot)\equiv s^*$
starting from
initial states  $\bar s(t_0)$ and $s^*$ are unique. (Former FSP is the unique solution of ODE \eqn{eq-deriv-33},
and the uniqueness of the latter is by (iv).) Therefore, 
by (i), $\bar s(t) \ge s^*$ for $t\ge t_0$. 
Then $s^{**} \ge s^*$. 
Furthermore, it is easy to see from the ODE \eqn{eq-deriv-33} that $s^{**} \ne s^*$
is impossible; otherwise, $\sum_j \bar x_{1,j}(t) \downarrow \sum_j x_{1,j}^{**} > \sum_j x_{1,j}^*$, which would imply that 
$(d/dt)\sum_j \bar x_{1,j}$ must stay negative bounded away from zero for all $t$.
We have proved the convergence $\bar s(t) \to s^*$.

To prove uniqueness of $\bar s(\cdot)$, fix any such FSP.
We will compare it to the FSP $s^{(\delta)}(\cdot)$, 
with the following initial state: $\xi_{r,j}^{(\delta)}(0) = \mu_j \bar y_{1,j}(0) \delta = \mu_j \beta_j \delta$,
where $\delta>0$ is a parameter. (Note that $\mu_j \bar y_{1,j}(0)$ is what we denoted $w_j(0)$ 
in the proof of (v); in other words, $\xi_{r,j}^{(\delta)}(0)$ are small numbers, proportional to $b_j$.)
For a fixed small $\delta$, the conditions of (iii) hold for
FSP $s^{(\delta)}(\cdot)$ and, therefore, it is unique in a sufficiently small interval $[0,\epsilon]$. Moreover, the choice
of $\xi_{r,j}^{(\delta)}(0)$ guarantees that $(d/dt) \xi_{r,j}^{(\delta)}(t) >0$ in that interval, which means 
that $s^{(\delta)}(\cdot)$ is monotone non-increasing in $[0,\epsilon]$. We can repeat the arguments, applied 
above to $\bar s(\cdot)$, to establish that FSP $s^{(\delta)}(\cdot)$ is unique monotone non-increasing
in $[0,\infty)$, and $s^{(\delta)}(t) \to s^*$. Further note that for any sufficiently small $t_0>0$,
$s^{(\delta)}(t_0) \le \bar s(t_0)$. By (i), $s^{(\delta)}(t) \le \bar s(t)$ for $t\ge t_0$, and then for 
all $t\ge 0$. Therefore, $\xi_{r,j}^{(\delta)}(t) \ge \bar \xi_{r,j}(t)$ for all $t$ and all $(r,j)$,
and the sum $h(t) = \sum_{r,j} (\xi_{r,j}^{(\delta)}(t) - \bar \xi_{r,j}(t))$ can serve as a ``distance'' 
between $s^{(\delta)}(t)$ and $\bar s(t)$. From ODE \eqn{eq-deriv-33} we obtain
$$
\frac{dh}{dt} = - \sum_j \mu_j \sum_r (\xi_{r,j}^{(\delta)}(t) - \bar \xi_{r,j}(t)) \le 0, ~~t > 0.
$$
so that $h(t)$ is non-increasing, and $h(t) \le h(0) = R \mu_j \beta_j \delta$.
This imples that $s^{(\delta)}(\cdot) \to \bar s(\cdot)$ u.o.c. as $\delta\downarrow 0$.
Since $s^{(\delta)}(\cdot)$ is unique for each $\delta$, this proves uniqueness of $\bar s(\cdot)$.

(vii) From (ii) and definition of $\nu_j$, using relation $y_{1,j}(t)=x_{1,j}(t)-x_{2,j}(t)$,
we have
$$
(d/dt) x_{1,j}(t) = \mu_j x_{2,j}(t), ~~j\in \cj.
$$
(For $t=0$ -- it is the right derivative.)
Also from (ii), we observe that 
in a sufficiently small fixed neighborhood of time $t$,
 the derivative $(d/du) x_{1,\ell}(u)$
is uniformly Lipschitz continuous (given the statement assumptions which hold at $t$).
This implies that, for an arbitrarily small $\epsilon_1>0$,
in a (further reduced) small neighborhood of $t$,
$(d/du) x_{1,\ell}(u) \ge \mu_j \epsilon -\epsilon_1$;
which in turn implies the desired property.
$\Box$

\section{System with unit buffer sizes}
\label{sec-unit-buffers}

Here we consider a special system where all buffer sizes are equal to $1$, i.e. $B_j=1$ for all $j$.
The process $s^n(\cdot)$ is automatically positive recurrent (Lemma~\ref{lem-stabil-finite}).

\begin{lem}
\label{lem-lower-bound11}
Suppose $B_j=1$ for all $j$. Then
\beql{eq-lem1111}
s^n(\infty) \Rightarrow s^*.
\end{equation}
\end{lem}

{\em Proof.} Since space $\cs$ is compact, any subsequence of $n$ has a further subsequence,
along which
\beql{eq-stationary-conv999}
s^n(\infty) \Rightarrow s^{\circ}(\infty),
\end{equation}
where $s^{\circ}(\infty)$ is a random element in $\cs$.
We need to prove that 
\beql{eq-lem11}
s^{\circ}(\infty) = s^*, ~~w.p.1.
\end{equation}
For each $n$, consider the process $s^n(\cdot)$, 
starting from the maximum initial state: 
$\xi^n_{r,j}=0$ for all $j$ and $r$.
By Corollary~\ref{cor-monotone}, process $s^n(\cdot)$ is monotone non-increasing.
Consider any fixed pair $(r,j)$.
Fix arbitrary $\epsilon>0$, and choose $T>0$ large enough so that
the FSP $\bar s(\cdot)$ starting from the maximum initial condition (as in
Lemma~\ref{lem-ms-fsp-prop}(vi)) is such that $\bar\xi_{r,j}(T) \ge \xi_{r,j}^* -\epsilon/2$.
Then, by Lemma~\ref{lem-ms-fsp-conv}, $\pr\{\xi_{r,j}^n(T) \ge \xi_{r,j}^* -\epsilon\} \to 1$.
From here and the fact that $s^n(\cdot)$ is monotone non-increasing,
we obtain
$$
\liminf_{n\to\infty} \pr\{\xi_{r,j}^n(\infty) \ge \xi_{r,j}^* -\epsilon\} \ge 
\liminf_{n\to\infty} \pr\{\xi_{r,j}^n(T) \ge \xi_{r,j}^* - \epsilon\} = 1.
$$
Therefore, since $\{\xi_{r,j} \ge \xi_{r,j}^* -\epsilon\}$ is a closed set, 
$$
\pr\{\xi_{r,j}^{\circ}(\infty) \ge \xi_{r,j}^* - \epsilon\} \ge 
\limsup_{n\to\infty} \pr\{\xi_{r,j}^n(\infty) \ge \xi_{r,j}^* -\epsilon\}
\ge 1.
$$
This holds for any $\epsilon>0$, so we have 
$\pr\{\xi_{r,j}^{\circ}(\infty) \ge \xi_{r,j}^* \} = 1$ for any pair $(r,j)$.
This is equivalent to (recall that here we consider system with $B_j=1$ for all $j$)
\beql{eq-temp-lower}
s^{\circ}(\infty) \le_{st} s^*. 
\end{equation}
Note that \eqn{eq-temp-lower} implies $x_{1,j}^{\circ}(\infty) \le x_{1,j}^*=\nu_j, ~\forall j$, w.p.1.
To prove lemma, it remains to show that, in fact, \eqn{eq-lem11} holds.

For a system with fixed parameter $n$,
the (fluid-scaled) steady-state rate at which customers enter service (equal to the 
rate at which they leave service) is lower bounded by
$$
\sum_r (\lambda/R) \pr\{\sum_j \xi_{r,j}^n(\infty) >0\},
$$
which by \eqn{eq-temp-lower} converges to $\lambda$ as $n\to\infty$.
But, this rate is obviously upper bounded by $\lambda$.
Therefore, for the steady-state rate at which customers leave service, we must have
$$
\lim_{n\to\infty} \E \sum_j \mu_j x_{1,j}^n(\infty) = \lambda.
$$
But,
$$
\lim_{n\to\infty} \E \sum_j \mu_j x_{1,j}^n(\infty) = \E \sum_j \mu_j x_{1,j}^{\circ}(\infty).
$$
We conclude that
\beql{eq-temp-rate-lam}
\E \sum_j \mu_j x_{1,j}^{\circ}(\infty) = \lambda.
\end{equation}

Suppose now that \eqn{eq-temp-lower} holds, but \eqn{eq-lem11} does not;
namely, for at least one $j$, we have 
$\pr\{x_{1,j}^{\circ}(\infty) < \nu_j\} >0$. 
This is impossible, because then
$
\E \sum_j \mu_j x_{1,j}^{\circ}(\infty) < \sum_j \mu_j \nu_j = \lambda
$, which contradicts \eqn{eq-temp-rate-lam}. This proves \eqn{eq-lem11}. $\Box$

\section{Proof of Theorem~\ref{thm-stabil-infinite}}
\label{sec-stability}

Recall that we consider a system with fixed buffer sizes $B_j$, some of which are infinite. 
Denote by $\cj_\infty$ the subset of pools $j$ with $B_j=\infty$.

\subsection{Instability characterization.}
\label{sec-instab-char}

Here we consider a system with the scaling parameter $n$ being fixed. System {\em instability} means the 
originally defined process (i.e., with original -- not compactified -- state space) is not ergodic
 (i.e., not positive recurrent).

Our goal is characterize instability, specifically to prove Lemma~\ref{lem-305} at the end of this subsection.
The statement of Lemma~\ref{lem-305} involves three processes: $S^n(\cdot)$, which is the main process,
whose instability we try to characterize; $\hat S^n(\cdot)$ is an auxiliary process, for a system such that all queues in one of the pools $j$ (with $B_j=\infty$) are infinite, while buffer sizes in all other pools $\ell \ne j$
are $1$; $\breve S^n(\cdot)$ is an auxiliary process (considered in Lemma~\ref{lem-lower-bound11}), 
for a system such that all buffer sizes in all pools
are $1$. Processes $s^n(\cdot)$, $\hat s^n(\cdot)$ and $\breve s^n(\cdot)$ are the corresponding
mean-field versions.
Informally speaking, Lemma~\ref{lem-305} states that if the process $s^n(\cdot)$ is unstable, then there exists
pool $j$ such that, as $t\to\infty$, this pool becomes fully occupied, and furthermore, all queues in that pool
converge to infinity; moreover, ``in the $t\to\infty$ limit'', $s^n(\infty)$ dominates $\hat s^n(\infty)$, which in turn 
dominates $\breve s^n(\infty)$. Lemmas~\ref{lem-301}-\ref{lem-303} sequentially strengthen one another,
to show (in Lemmas~\ref{lem-303}) that instability implies that, in at least one pool $j$, all queues must 
converge to infinity (in probability) as $t\to\infty$. The  key Lemma~\ref{lem-304} uses this fact to show that,
``in the $t\to\infty$ limit'', $S^n(\infty)$ dominates $\hat S^n(\infty)$; 
to do that, in the proof of Lemma~\ref{lem-304} yet another auxiliary process, $\tilde S^n(\cdot)$,
is used, which coincides with $\tilde S^n(\cdot)$ up to some time point, but then switches to a ``smaller'' process
with buffer sizes $B_\ell=1$ in the pools $\ell \ne j$.
Finally, 
using the dominance $\breve S^n(\infty) \le_{st} \hat S^n(\infty)$ (which follows from monotonicity), and
switching to mean-field processes, we obtain  Lemma~\ref{lem-305} as a corollary of  Lemma~\ref{lem-304}.

\begin{lem}
\label{lem-301}
Suppose the system is unstable. Then there exists $j\in \cj_\infty$ and $\delta>0$ such that the following holds. For any $C>0$, there exists $T=T(C)>0$ such that, uniformly on all initial states,
for any $i\in\cn_j$,
\beql{eq-301}
\pr\{Q_i^n(T) \ge C\} \ge \delta.
\end{equation}
\end{lem}

{\em Proof.} It suffices to prove that the lemma statement holds uniformly on all {\em idle} initial states (i.e., such that $Q_i^n(0), ~\forall i$). This follows from the fact that any initial state dominates {\em one of the
idle states}, and from Lemma~\ref{lem-monotone}. 
Pick arbitrary proper state  $\sigma^n_0$ (with all queues finite) of the process $S^n(\cdot)$.
There is only a finite number of idle states, and the Markov process $S^n(\cdot)$ is irreducible.
Then,
for arbitrary fixed $T_0>0$, uniformly on the idle initial states, 
with probability at least some $\delta_0>0$,
the system in time $T_0$ makes transition to state $\sigma^n_0$.
Therefore, it suffices to prove the lemma statement for the single initial state $S^n(0)=\sigma^n_0$;
this is what we can and do assume in the rest of this proof.

The system instability by definition means that the (irreducible) process $S^n(\cdot)$ is not positive recurrent
(and has no stationary distribution). Then, as $t\to\infty$, the probability of being in any
given state goes to $0$ (cf. \cite{Grim_Stir}, Theorem (6.9.21), p. 261).
Therefore,  $\sum_i Q_i^n(T) \stackrel{P}{\rightarrow} \infty$ as $T\to\infty$
(because for any constant $K$ there is only a finite number of states with $\sum_i Q_i^n \le K$).
Then, for any $C>0$, we can pick a sufficiently large $T>0$, such that 
$$
\pr\{\sum_i Q_i^n(T) \ge Cn\} \ge 1/2.
$$
But,
$$
\pr\{\sum_i Q_i^n(T) \ge Cn\} \le \sum_i \pr\{Q_i^n(T) \ge C\}.
$$
Then, there exists $j$ and $i\in \cn_j$ for which $\pr\{Q_i^n(T) \ge C\} \ge \delta=1/(2n)$;
then, by symmetry, this is true for all $i\in \cn_j$.
This is true for a sequence of $C\uparrow\infty$ and the corresponding sequence of $T=T(C)$;
we choose a subsequence, for which the corresponding $j$ is same. Clearly, this $j\in \cj_\infty$.
Therefore, the desired property holds for this $j$ and $\delta=1/(2n)$.
$\Box$

\begin{lem}
\label{lem-302}
Suppose the system is unstable. Consider $j\in \cj_\infty$ picked as in Lemma~\ref{lem-301}.
Then for some $\delta>0$ and for any $C>0$ there exists $T=T(C)>0$ such that, uniformly on all initial states,
\beql{eq-302}
\pr\{\min_{i\in \cn_j} Q_i^n(T) \ge C\} \ge \delta.
\end{equation}
\end{lem}

{\em Proof.} The proof is by induction on subsets of $\cn_j$. Suppose the lemma statement holds for 
subsets $\hat \cn_j \subset \cn_j$ of cardinality at most $\kappa < |\cn_j|$. Namely,  for some $\delta>0$ and 
for any $C>0$ there exists $T=T(C)>0$ such that, uniformly on all initial states, for any 
subset $\hat \cn_j \subset \cn_j$ with $|\hat \cn_j| \le \kappa$,
\beql{eq-302a}
\pr\{\min_{i\in \hat \cn_j} Q_i^n(T) \ge C\} \ge \delta.
\end{equation}
Note that by Lemma~\ref{lem-301} this is true for $\kappa=1$; this is the induction base. 
Pick any $\hat \cn_j \subset \cn_j$ with $|\hat \cn_j| = \kappa$, and an element
$m\in \cn_j \setminus \hat \cn_j$. Since there is only a finite number of subsets of $\cn_j$,
the induction step will be completed if we show that 
for some $\delta'>0$ and 
for any $C>0$ there exists $T'=T'(C)>0$ such that, uniformly on all initial states,
\beql{eq-302b}
\pr\{\min_{i\in \hat \cn_j \cup \{m\}} Q_i^n(T') \ge C \} \ge \delta'.
\end{equation}
Pick $\delta$, arbitrarily large $C$ and the corresponding $T=T(C)$, so that 
\eqn{eq-302a} holds
for the subset $\hat \cn_j$. Pick any $\epsilon \in (0,\delta)$. For the same $\delta$, 
we can pick $C_1>0$ so large that, if $Q_m^n(0)\ge C_1$ then
$$
\pr\{Q_m^n(T) \ge C \} \ge 1-\epsilon.
$$
We can and do pick $T_1=T_1(C_1)$, so that \eqn{eq-301} holds with $i=m$, $T=T_1$, and $C=C_1$.

With arbitrary fixed initial state, we see that 
$\pr\{Q_m^n(T_1) \ge C_1 \} \ge \delta$. 
We also have
$$
\pr\{\min_{i\in \hat \cn_j} Q_i^n(T_1+T) \ge C ~|~ Q_m^n(T_1) \ge C_1\} \ge \delta
$$
and
$$
\pr\{Q_m^n(T_1+T) \ge C ~|~ Q_m^n(T_1) \ge C_1\} \ge 1-\epsilon.
$$
Therefore,
$$
\pr\{\min_{i\in \hat \cn_j \cup \{m\}} Q_i^n(T_1+T) \ge C \} \ge \delta (\delta - \epsilon).
$$
This proves the induction step \eqn{eq-302b}, with $\delta'=\delta (\delta - \epsilon)$, and $T'=T_1+T$ chosen, as a function 
of $C$, as described above. $\Box$

\begin{lem}
\label{lem-303}
Suppose the system is unstable. Consider $j\in \cj_\infty$ picked as in Lemma~\ref{lem-301}
(and Lemma~\ref{lem-302}).
Then, uniformly on all initial states,
\beql{eq-303}
\min_{i\in \cn_j} Q_i^n(t) \stackrel{P}{\rightarrow} \infty, ~\mbox{as}~ t\to\infty.
\end{equation}
\end{lem}

{\em Proof.} 
First, we show that if the statement of Lemma~\ref{lem-302} holds for some fixed $\delta>0$, then it must
hold for any $\delta \in (0,1)$. To show this, it suffices to verify the following fact: 
if the statement of Lemma~\ref{lem-302} holds for some fixed $\delta>0$, then it also holds
for $\delta'=\delta(1-\epsilon) + (1-\delta) \delta=\delta(2-\delta-\epsilon)$ with arbitrarily small fixed $\epsilon>0$. We will do just that.

Fix arbitrary initial state. Fix $\epsilon>0$. Pick $\delta$ as in Lemma~\ref{lem-302}, 
then pick arbitrarily large $C$ and the corresponding $T=T(C)$, so that \eqn{eq-302} holds.
Choose $C_1>0$ so large that, if
$\min_{i\in \cn_j} Q_i^n(0) \ge C_1$ then
$$
\pr\{\min_{i\in \cn_j} Q_i^n(T) \ge C\} \ge 1-\epsilon.
$$
Then, again using Lemma~\ref{lem-302}, choose $T_1=T_1(C_1)$ so that \eqn{eq-302} holds 
with $C$ and $T$ replaced by $C_1$ and $T_1$, respectively.
We have
$$
\pr\{\min_{i\in \cn_j} Q_i^n(T_1) \ge C_1\} \ge \delta.
$$
Considering the system state at time $T_1+T$, conditioning on the event $\{\min_{i\in \cn_j} Q_i^n(T_1) \ge C_1\}$ occurring or not, we obtain
$$
\pr\{\min_{i\in \cn_j} Q_i^n(T_1+T) \ge C\} \ge \delta(1-\epsilon) + (1-\delta) \delta = \delta'.
$$
Thus, the statement of Lemma~\ref{lem-302} indeed holds with $\delta$ rechosen to be $\delta'$
and function $T(C)$ rechosen to be $T'(C)$. Consequently, the statement of Lemma~\ref{lem-302} holds 
for any $\delta \in (0,1)$.

Then, the statement of this lemma easily follows, because for any $\delta\in (0,1)$
and arbitrarily large $C>0$ we can choose $T=T(C)$ such that, for any initial state
$$
\pr\{\min_{i\in \cn_j} Q_i^n(t) \ge C\} \ge \delta
$$
for $t=T$. But then this is also true for all $t\ge T$.
$\Box$

Consider $j\in \cj_\infty$ picked as in Lemma~\ref{lem-301}
(and Lemmas~\ref{lem-302} and \ref{lem-303}).
Consider a variant of the system, such that all queues in pool $j$ are infinite, $Q_i^n(t)\equiv \infty$,
$i \in \cn_j$, and buffer sizes in all other pools are equal to $1$, $B_\ell = 1$, $\ell \ne j$.
Recall that if a queue size is infinite at time $0$, it remains infinite forever.
Therefore, the process for this specific system, let us denote it $\hat S^n(\cdot)$, is a Markov
chain with a {\em finite} state space.
It has unique stationary distribution; then $\hat S^n(\infty)$ is its random state in stationary regime.
Since the state space is finite, the convergence to the stationary distribution is 
uniform on the initial states:
$$
\lim_{t\to\infty} \max_{\hat S^n(0)} 
\|Dist[\hat S^n(t)] - Dist[\hat S^n(\infty)]\| = 0.
$$

\begin{lem}
\label{lem-304}
Suppose the system is unstable. Consider $j\in \cj_\infty$ picked as in Lemma~\ref{lem-301}
(and Lemmas~\ref{lem-302} and \ref{lem-303}). 
Let an arbitrary initial state $S^n(0)$ be fixed.
Consider a random element $S^n(\infty)$, which is any limit in distribution of $S^n(t)$ along a subsequence of $t\to\infty$.
(Recall that our process is viewed as having the compact state space, with possibly infinite queues.
Then, any subsequence of $t\to\infty$ has a further subsequence, along which  $S^n(t)$ converges in distribution.)
Then, $\hat S^n(\infty) \le_{st} S^n(\infty)$. (Note that this implies $Q^n_i(\infty) = \infty, ~\forall i \in \cn_j$.)
\end{lem}

{\em Proof.} 
It will suffice to prove the following property; let us label it (P1):\\
{\em For any $\epsilon_1>0$, $\epsilon_2>0$ and any $C>0$,
there exists $T'=T'(C)$ 
such that the following holds uniformly on all initial states $S^n(0)$.
Process $S^n(\cdot)$ can be constructed on a common probability space with another process
$\tilde S^n(\cdot)$, so that
there exists an event $G$ with $\pr\{G\} \ge 1-\epsilon_1$, 
such that
$$
\pr\{ \min_{i\in \cn_j}  Q_i^n(T') \ge C  ~|~ G\} =1,
$$
$$
\tilde S^n(T') ~\mbox{is dominated by}~ S^n(T') ~\mbox{on the subset of servers $\cn \setminus \cn_j$},
~~\mbox{on the event $G$},
$$
$$
\|Dist[\tilde Q^n(T')] - Dist[\hat Q^n(\infty)]  \| \le \epsilon_2.
$$
}
We proceed with proving property (P1).
Let us first choose $T_2>0$ sufficiently large so that, uniformly on the initial states $\hat S^n(0)$ of the process 
$\hat S^n(\cdot)$, 
\beql{eq-conv-hat}
\|Dist[\hat Q^n(T_2) - Dist[\hat Q^n(\infty)]  \| \le \epsilon_2.
\end{equation}
Then, we choose $C_1>C$, sufficiently large so that if $\min_{i\in \cn_j} Q^n_i(0) \ge C_1$ then
with probability at least $1-\epsilon_1/2$, $\min_{i\in \cn_j} \min_{0\le t \le T_2} Q^n_i(t) \ge C$. Then we choose $T_1=T_1(C_1)>0$
sufficiently large so that, uniformly on the initial states $S^n(0)$, the event
$$
\{ \min_{i \in \cn_j} Q^n_i(T_1) \ge C_1 \}
$$
occurs with probability at least $1-\epsilon_1/2$. We choose $T'=T_1+T_2$ and event $G$ to be
$$
G=\{ \min_{i\in \cn_j} Q^n_i(T_1) \ge C_1, ~~ \min_{i\in \cn_j} \min_{T_1 \le t \le T'} Q^n_i(t) \ge C\}.
$$
We construct the process $\tilde S^n(\cdot)$ as follows. In the time interval $[0,T_1)$, 
$\tilde S^n(\cdot)$ coincides with $S^n(\cdot)$. At time $T_1$ we change the model for the process
$\tilde S^n(\cdot)$ by keeping buffer sizes in pool $j$ to be infinite, $\tilde B_j = \infty$, and making buffer sizes
in all other pools to be equal to $1$, $\tilde B_\ell = 1$ for $\ell \ne j$; we also change all queue sizes in 
pool $j$ to be infinite, $\tilde Q^n_i(T_1)=\infty, ~i\in\cn_j$, and leaving at most one customer in
each queue in the other pools, 
$$
\tilde Q^n_i(T_1)=1 \wedge \tilde Q^n_i(T_1-), ~i\in\cn_\ell, ~ \ell \ne j.
$$
Starting time $T_1$, the process $\tilde S^n(\cdot)$ is coupled to the process $S^n(\cdot)$,
as described in Lemma~\ref{lem-monotone},
until the random time $\tau$, $T_1 \le \tau \le \infty$, when $\min_{i\in \cn_j} Q^n_i(t) = 0$ for the first time 
after $T_1$. Obviously, $\tau > T'$ on the event $G$. It remains to observe that
in the interval $[T_1,\tau)$, the process $S^n(\cdot)$ dominates $\tilde S^n(\cdot)$ 
on the subset $\cn \setminus \cn_j$ of servers.
$\Box$

Observe that $\breve S^n(\infty) \le_{st} \hat S^n(\infty)$, where 
$\breve S^n(\cdot)$ is the
 process for the system with $B_\ell=1$ for all $\ell$, considered in Lemma~\ref{lem-lower-bound11}.
Let $\breve s^n(\infty)$ and  $\hat s^n(\infty)$ be the random states (in stationary regime)
for the corresponding mean-field (or, fluid-scaled) processes. 
Also note that from the definition of process $\hat s^n(\cdot)$, 
$\hat x_{1,j}^n(\infty) = \beta_j$ for all $n$.
Then, as a corollary of Lemma~\ref{lem-304}
we obtain the following 

\begin{lem}
\label{lem-305}
Suppose the system is unstable. Consider $j\in \cj_\infty$ picked as in Lemma~\ref{lem-301}
(and Lemmas~\ref{lem-302} -- \ref{lem-304}). 
Let an arbitrary initial state $s^n(0)$ be fixed.
Consider a random element $s^n(\infty)$, which is any limit in distribution of $s^n(t)$ along a subsequence of $t\to\infty$.
Then, $\breve s^n(\infty) \le_{st} \hat s^n(\infty) \le_{st} s^n(\infty)$ and 
$x_{1,j}^n(\infty) = \hat x_{1,j}^n(\infty) = \beta_j$.
\end{lem}

\subsection{Theorem~\ref{thm-stabil-infinite} proof.}
\label{sec-proof-stability}

For each $n$, consider the process, starting from an arbitrarily chosen
 fixed initial state $s^n(0)$ (with all queues being finite). 
We have
\beql{eq-temp555}
\limsup_{T\to\infty} (1/T) \int_0^T [\E \sum_\ell \mu_\ell x_{1,\ell}^n(t)]dt \le \lambda,
\end{equation}
because the LHS is the $\limsup$ of the time-averaged expected (fluid-scaled) number of customer
service completions over the interval $[0,T]$, which cannot exceed the limit of the time-averaged expected (fluid-scaled) number of customer arrivals.

Suppose the system is unstable 
for infinitely many $n$. 
This leads to a contradiction as follows.
Consider an infinite subsequence of those $n$, for which the system is unstable.
Moreover, we choose, if necessary, a further subsequence, along which 
we can apply Lemma~\ref{lem-305} with the same value of $j$.
For each such $n$, it follows from Lemma~\ref{lem-305} that
\beql{eq-aaa1}
\lim_{t\to\infty} \E x_{1,j}^n(t) = \beta_j,
\end{equation}
\beql{eq-aaa2}
\liminf_{t\to\infty} \E x_{1,\ell}^n(t) \ge \E \breve x_{1,\ell}^n(\infty), ~~\ell \ne j.
\end{equation}
(Indeed,  consider, say, \eqn{eq-aaa2} for a fixed $\ell \ne j$. 
By Lemma~\ref{lem-305}, any sequence of $t$ increasing to infinity,
contains a further subsequence
along which $s^n(t) \Rightarrow s^n(\infty)$ and \eqn{eq-aaa2} holds.)
Then,
\beql{eq-temp666}
\liminf_{T\to\infty} (1/T) \int_0^T [\E \sum_\ell \mu_j x_{1,\ell}^n(t)]dt \ge 
\mu_j \beta_j +\sum_{\ell \ne j} \mu_\ell \E  \breve x_{1,\ell}^n(\infty).
\end{equation}
Recall that, by Lemma~\ref{lem-lower-bound11}, $\breve s^n(\infty) \Rightarrow s^*$
as $n\to\infty$. This implies that, as $n\to\infty$, the RHS in \eqn{eq-temp666}
converges to $\mu_j \beta_j + \sum_{\ell \ne j}\mu_\ell \nu_\ell > \lambda$. This contradicts the fact
that \eqn{eq-temp555} must hold for each $n$.
$\Box$

\section{Proof of Theorem~\ref{thm-infinite}}
\label{sec-proof-finite-buffers}

Space $\cs$ is compact. Therefore, any subsequence of $n$ has a further subsequence,
along which
\beql{eq-stationary-conv}
s^n(\infty) \Rightarrow s^{\circ}(\infty),
\end{equation}
where $s^{\circ}(\infty)$ is a random element in $\cs$.
Therefore, to prove Theorem~\ref{thm-infinite} it suffices to show that
any limit in \eqn{eq-stationary-conv} is equal to $s^*$ w.p.1.

From Lemma~\ref{lem-lower-bound11}, 
we obtain the following Corollary~\ref{lem-lower-bound}, 
because $s^n(\infty)$ for a system with arbitrary buffer sizes $B_j \ge 1$ stochastically dominates that for the system with $B_j=1$ for all $j$, and therefore this relation holds also for the limit $s^{\circ}(\infty)$.

\begin{cor}
\label{lem-lower-bound}
Any subsequential limit $s^{\circ}(\infty)$ in \eqn {eq-stationary-conv}
is such that
$$
s^* \le s^{\circ}(\infty), ~~w.p.1.
$$
\end{cor}

{\em Proof of Theorem~\ref{thm-infinite}.} 
Consider any subsequential limit $s^{\circ}(\infty)$  in \eqn{eq-stationary-conv},
along a subsequence of $n$; for  the rest of the proof, we consider this subsequence.
By Corollary~\ref{lem-lower-bound},
$$
\E \sum_j \mu_j x_{1,j}^{\circ}(\infty) \ge \lambda.
$$
On the other hand, for any $n$, the (fluid-scaled) steady-state rate of service completions in the system is 
$$
\E \sum_j \mu_j x_{1,j}^n(\infty) \le \lambda;
$$
this implies
$$
\E \sum_j \mu_j x_{1,j}^{\circ}(\infty) = \lim_{n\to\infty} \E \sum_j \mu_j x_{1,j}^n(\infty) \le \lambda.
$$
Therefore,
$$
\E \sum_j \mu_j x_{1,j}^{\circ}(\infty) = \lambda,
$$
which (by Corollary~\ref{lem-lower-bound}) implies
\beql{eq-x1-is-nu}
x_{1,j}^{\circ}(\infty) = \nu_j, ~j \in \cj, ~~~w.p.1.
\end{equation}
Now, \eqn{eq-x1-is-nu} means that, w.p.1 for each $j$, $\sum_r \xi_{r,j}^{\circ}(\infty) = \sum_r \xi_{r,j}^*$,
which along with Corollary~\ref{lem-lower-bound} implies
$$
\xi_{r,j}^{\circ}(\infty) = \xi_{r,j}^*, ~~j \in \cj,~ r\in\calr, ~~~w.p.1.
$$
It remains to show that 
\beql{eq-x2-is-0}
x_{2,j}^{\circ}(\infty) = 0, ~j \in \cj, ~~~w.p.1.
\end{equation}
Suppose not, namely, for at least one $\ell$, 
$\pr\{x_{2,\ell}^{\circ}(\infty) > \epsilon\} = 2\epsilon_1$, for some 
$\epsilon>0$, $\epsilon_1>0$.
Then, for all sufficiently large $n$, 
$\pr\{x_{2,\ell}^n(\infty) > \epsilon\} > \epsilon_1$. 
For each sufficiently large $n$, consider a stationary version of $s^n(\cdot)$,
that is $s^n(t)\stackrel{d}{=}s^n(\infty)$ for all $t\ge 0$.
Then $\pr\{x_{2,\ell}^n(0) > \epsilon\} > \epsilon_1$. Now, employing 
Lemma~\ref{lem-ms-fsp-conv} and Lemma~\ref{lem-ms-fsp-prop}(vii), we can easily
show that, for some $\tau>0$, $\delta>0$, and all large $n$,
$$
\pr\{x_{1,\ell}^n(\tau) \ge \nu_\ell + \delta/2\} > \epsilon_1/2.
$$
But then 
$$
\pr\{x_{1,\ell}^{\circ}(\infty) \ge \nu_\ell+ \delta/2\} \ge \limsup_{n\to\infty}
\pr\{x_{1,\ell}^n(\tau) \ge \nu_\ell + \delta/2\} \ge \epsilon_1/2,
$$
a contradiction with \eqn{eq-x1-is-nu}, which proves \eqn{eq-x2-is-0}.
$\Box$

\section{Conclusions}
\label{sec-conclusion}

The main conclusion of this paper is that a pull-based approach for load distribution can be applied in 
large-scale heterogeneous systems with multiple independent routers. 
We proposed (two versions of) a specific algorithm PULL -- along with a specific mechanism for its implementation --
and proved its asymptotic optimality. PULL algorithm in this paper is a generalization 
of the single-router  algorithm in \cite{St2014_pull}. The generalization of the algorithm itself
is very natural, but a priori it is far from obvious that it is provably optimal for the multi-router model.
In fact, to make our proofs work, the version PULL-2 of the algorithm is such that,
in addition to pull-messages, it requires occasional ``pull-remove-messages'' to keep the pull-messages 
at the routers ``up to date.'' This additional mechanism potentially increases the router-server message exchange
rate to two-per-customer; however, as our results show, in the large-scale asymptotic limit, pull-remove-messages are never used,
and therefore the limiting exchange rate is only one-per-customer, same as in the single-router case.
We conjecture that the additional pull-remove-message mechanism is {\em not} necessary for
PULL-2 asymptotic optimality to hold. However, its absence substantially complicates the analysis.
Verifying this conjecture may be one of the subjects of future work.

\bibliographystyle{acmtrans-ims}
\bibliography{biblio-stolyar}

\begin{thebibliography}{}
\ifx \url   \undefined \def \url#1{#1}   \fi

\bibitem{BB08}
\textsc{Badonnel, R.} \textsc{and} \textsc{Burgess, M.} (2008).
\newblock Dynamic pull-based load balancing for autonomic servers.
\newblock \emph{Network Operations and Management Symposium, NOMS 2008\/},
  751--754.

\bibitem{BLP2012-jsq-asymp-indep}
\textsc{Bramson, M.}, \textsc{Lu, Y.}, \textsc{and} \textsc{Prabhakar, B.}
  (2012).
\newblock Asymptotic independence of queues under randomized load balancing.
\newblock \emph{Queueing Systems\/}~\emph{71}, 247--292.

\bibitem{BLP2013-jsq-asymp-tail}
\textsc{Bramson, M.}, \textsc{Lu, Y.}, \textsc{and} \textsc{Prabhakar, B.}
  (2013).
\newblock Decay of tails at equlibrium for fifo join the shortest queue
  networks.
\newblock \emph{The Annals of Applied Probability\/}~\emph{23}, 1841--1878.

\bibitem{EsGam2015}
\textsc{Eschenfeldt, P.} \textsc{and} \textsc{Gamarnik, D.} (2015).
\newblock Join the shortest queue with many servers. the heavy traffic
  asymptotics.

\bibitem{Grim_Stir}
\textsc{Grimmett, G.} \textsc{and} \textsc{Stirzaker, D.} (2001).
\newblock \emph{Probability and Random Processes (3rd ed.)}.
\newblock Oxford University Press.

\bibitem{Liggett-book}
\textsc{Liggett, T.~M.} (1985).
\newblock \emph{Interacting Particle Systems}.
\newblock Springer.

\bibitem{G11}
\textsc{Lu, Y.}, \textsc{Xie, Q.}, \textsc{Kliot, G.}, \textsc{Geller, A.},
  \textsc{Larus, J.}, \textsc{and} \textsc{Greenberg, A.} (2011).
\newblock Join-idle-queue: A novel load balancing algorithm for dynamically
  scalable web services.
\newblock \emph{Performance Evaluation\/}~\emph{68}, 1057--1071.

\bibitem{Mitz2001}
\textsc{Mitzenmacher, M.} (2001).
\newblock The power of two choices in randomized load balancing.
\newblock \emph{IEEE Transactions on Parallel and Distributed
  Systems\/}~\textbf{12},~10, 1094--1104.

\bibitem{MBLW2015}
\textsc{Mukherjee, D.}, \textsc{Borst, S.}, \textsc{van Leeuwaarden, J.},
  \textsc{and} \textsc{Whiting, P.} (2015).
\newblock Universality of load balancing schemes on diffusion scale.

\bibitem{St2015_grand-het}
\textsc{Stolyar, A.~L.} (2015a).
\newblock Large-scale heterogeneous service systems with general packing
  constraints.
\newblock arXiv:1508.07512.

\bibitem{St2014_pull}
\textsc{Stolyar, A.~L.} (2015b).
\newblock Pull-based load distribution in large-scale heterogeneous service
  systems.
\newblock \emph{Queueing Systems\/}~\textbf{80},~4, 341--361.

\bibitem{VDK96}
\textsc{Vvedenskaya, N.}, \textsc{Dobrushin, R.}, \textsc{and}
  \textsc{Karpelevich, F.} (1996).
\newblock Queueing system with selection of the shortest of two queues: an
  asymptotic approach.
\newblock \emph{Problems of Information Transmission\/}~\textbf{32},~1, 20--34.

\end{thebibliography}


\end{document}